\documentclass[preprint,longtitle,12pt,3p]{elsarticle}

\usepackage[utf8]{inputenc}

\usepackage{amssymb}
\usepackage{mathtools}
\usepackage[makeroom]{cancel}
\usepackage{amsthm}

\usepackage{amsfonts,amsmath}
\usepackage{graphicx}
\usepackage{algorithmicx}
\usepackage{algorithm}
\usepackage{algcompatible}
\usepackage{mathdots}
\usepackage{multirow}
\usepackage{accents}
\usepackage{siunitx}
\usepackage{dsfont}
\usepackage{multirow}
\newcommand{\Real}{\mathds{R}}

\usepackage{pgfplots}
\usepackage{pgfplotstable}
\usepackage{booktabs}
\usepackage{longtable}
\usepackage{tikz}
\usepackage{subfig}
\usetikzlibrary{arrows,shapes,mindmap,trees}
\usepackage[fancy,color=blue!70]{tikz-inet}
\usepackage{tikz-timing}
\usepackage{filecontents}
\usepackage{bm}
\usepackage{url}

\usepackage{array}

\usepackage{threeparttable}

\pgfplotsset{compat=1.6}

\newdefinition{remark}{Remark}

\journal{Electric Power Systems Research. \copyright 2020 CC-BY-NC-ND 4.0}

\begin{document}

\begin{frontmatter}



\title{Benders' decomposition of the unit commitment problem with semidefinite relaxation of AC power flow constraints}


\author[IBM]{M.~Paredes}
\ead{mparedes@br.ibm.com}

\author[IBM]{L.~S.~A.~Martins}
\ead{leonardo.martins@br.ibm.com}

\author[UNICAMP]{S.~Soares}
\ead{dino@cose.fee.unicamp.br}

\author[CSU]{Hongxing Ye}
\ead{h.ye@csuohio.edu}

\address[IBM]{IBM Research. São Paulo, Brazil.}
\address[UNICAMP]{University of Campinas. Campinas, Brazil.}
\address[CSU]{Cleveland State University. Cleveland, USA.}
\begin{abstract}
  In this paper we present a formulation of the unit commitment problem with AC power flow constraints. It is solved by a Benders’ decomposition in which the unit commitment master problem is formulated as a mixed-integer problem with linearization of the power generation constraints for improved convergence. Semidefinite programming relaxation of the rectangular AC optimal power flow is used in the subproblem, providing somewhat conservative cuts. Numerical case studies, including a 6-bus and the IEEE 118-bus network, are provided to test the effectiveness of our proposal. We show in our numerical experiments that the use of such strategy improves the quality of feasibility and optimality cuts generated by the solution of the convex relaxation of the subproblem, therefore reducing the number of iterations required for algorithm convergence.
\end{abstract}


\end{frontmatter}


\section{Introduction}%
\label{sec:introduction}

Unit commitment (UC) plays a central role in the operation of electric power systems and markets as the problem of finding the most economical on/off status and output dispatch of power generators connected to the network over the course of (most commonly) the next day~\cite{wood1984power}, such that physical and engineering constraints associated with power generation, transmission, and consumption components of the system are observed. It is a combinatorial optimization problem whose computational complexity scales exponentially with the size of the system under consideration. Moreover, nonlinearity of power flow physics makes this a challenging problem despite its long history in the academic literature~\cite{317549}, more so in light of the increasing trend in the penetration of non-dispatchable renewable energy sources (RES).

The taxonomy of unit commitment is broadly described in terms of its optimization constraints and solution methods, regardless of the context of its application, which usually defines whether the objective function reflects the economic welfare or costs associated with fuel, and unit startup and shutdown. With respect to its solution, given the inherent non-convex nature of the problem, algorithms proposed in the literature can be roughly classified, as cleverly noted by~\citet{Fattahi2017}, into ``methods based on a single convex [approximation] model, \ldots a series of convex models, [or] heuristics and local-search \ldots'' in which the key driver to their choice should be the set of constraints considered in the formulation. The basic set will have included unit ramping, minimum up- and down-time, limits on generation, and power balance equations in discrete time. Other constraint sets may also be formulated, such as reserve, security and, more generally, power flow constraints in either DC or (rectangular) AC formulations. More sophisticated models may also consider the stochasticity of RES.

\subsection{Related work}%
\label{sec:introduction:bibliography}

A properly comprehensive bibliographical review of solution methods proposed for the UC problem and its many variants is out of our intended scope for the paper. Therefore, this section traces back past work related to a lineage of algorithms of interest to AC transmission-constrained problem formulations, especially those employing convex relaxation. Since its inception, several techniques to solve the UC problem have been reported~\cite{1295033}, with dynamic programming (DP) having been explored in early works. \citet{4074522} presented a fully linearized formulation framework~\cite{4074522}, for which a so-called security function is proposed such that cold starts are penalized. This function has also been used to assess system security in an hourly, probabilistic basis~\cite{4074471}. \citet{4335130} and \citet{Kumar2007917} have considered ``\emph{a priori} information'' for the removal of infeasible paths to allow for computational tractability. Because of its combinatorial nature, the use of heuristics has then been shown to be suitable for large instances of the problem~\cite{5684940}, e.g. \citet{900083} proposed the use of greedy algorithms ordered by the average operation cost--these solutions are then evaluated for full supply of demand, and tested for feasibility by means of an optimal power flow (OPF).

In Lagrangian relaxation (LR) applications to UC it was common to relax reserve and demand coupling constraints in order to create individual subproblems for every unit. Other common approaches include the representation of individual unit subproblems as mixed-integer problems~\cite{Muckstadt1977} with techniques to select the Lagrange multipliers that maximize lower bounds produced by the relaxation, as well as the identification of identical solutions~\cite{1137620} that cause dual solutions to be far from the optimum by means of successive subproblems.

Other efficient methods to solve linear UC formulations were based on decomposition, such as Benders'~\cite{1405865}. One common approach of Benders' decomposition-based methods is to use piecewise-linear approximations of quadratic costs~\cite{Viana2013997}. \citet{1490608}, in one of the seminal works presenting AC transmission constraints, proposed to approximate the solution to the security-constrained UC master problem using LR and DP, with transmission constraints being checked in the subproblem in a rectangular AC power flow model, resulting in an optimization-simulation decomposition. Constraint violations are added back to the master problem as feasibility cuts in addition to cuts representing network contingencies. \citet{7061968} considered the uncertainty of wind power generation in real-time operation as the second of a two-stage stochastic program, with DC network-constrained day-ahead market clearing formulated in the first stage. In their proposed Benders' decomposition, the UC master problem is followed by a ``sufficiently  convexified'' AC network-constrained subproblem whose ``asymptotic convexification'' depends on the number of scenarios. \citet{Castillo2016} proposed a constraint-set linearization-based solution to the mixed-integer non-linear (MINLP) formulation of the UC problem with AC power flow constraints based on outer approximation (later extended by \citet{Liu2019} in a decomposition framework), where an analysis of the economic and operation advantages of AC network-constrained UC over its DC counterpart is provided.

Conic relaxation of polynomial optimization problems has driven interest in different applications in power systems, whether in the form of second-order conic (SOCP) or semidefinite (SDP) programs. Quadratic relaxations of binary variables~\cite{787402} and rectangular formulations of AC power flow~\cite{Molzahn2014,Bose2015,Huang2017} can be formulated as non-convex quadratically-constrained quadratic problems (QCQP) suitable for Shor's scheme-based SDP relaxations of MINLP formulations of the UC problem~\cite{4784461} with local resolution of relaxed infeasible solutions by means of approximation heuristics. \citet{Jabr201313} proposes a rank-1 constraint imposed on the SDP relaxation with an iterative convex rank reduction procedure. Local resolution of the relaxation is done by a greedy heuristic based on average fuel costs at full load in order to yield feasible schedules. \citet{Mhanna201285} propose the use of a heuristic based on a pruning algorithm, with complementary repair mechanisms that exploit the nature of problem constraints, and the SDP formulation itself. Another single convex model worth citing is that proposed in the work of \citet{Fattahi2017} for the DC network-constrained UC problem, in which reformulation-linearization is used to introduce relaxed non-convex quadratic inequalities in order to tighten the SDP relaxation seeking to avoid local feasibility resolution. More recently, \citet{Zohrizadeh2018} propose a power loss penalization method based on a series of SOCP relaxations that require an initialization ``sufficiently close to the feasible set'' for convergence.

\subsection{Contributions}%
\label{sec:introduction:contribution}

This paper contributes to the long line of methods for the solution of the UC problem building upon the recent developments in AC optimal power flow, and conic relaxation, with a QCQP formulation of the transmission constraints projected onto the space of semidefinite cones, for which the discrete unit commitment decisions are solved in a Benders' decomposition approach in order to overcome solutions with integrality gap due to a relaxation, while taking advantage of a series of convex SDP relaxations of the network. Furthermore, and in specific terms, the following contributions can be highlighted:
\begin{itemize}
	\item Not any sort of initialization with $\epsilon$-feasibility is required.
	\item A linearization of the power generation constraints is included in the master problem as a means to approximately represent the optimal power flow subproblem, and thus further constrain the search for integer solutions.
	\item In the proposed formulation, the rank of the SDP matrices associated with active power are guaranteed to be 1 at the end of every iteration of the algorithm.
	\item Consequently, a rank reduction procedure (RRP) is iteratively performed upon voltage variables, the only set of matrix variables for which rank-1 is not guaranteed.
\end{itemize}

Finally, two numerical case studies based on commonly available data sets are presented as evidence to illustrate the effectiveness of the proposed methodology.


\section{Unit Commitment Problem Formulation}%
\label{sec:problem}
We formulate the UC problem as a mixed-integer SDP problem of finding the unit commitment and active power generation dispatches that minimize costs subject to AC power flow constraints:
\begin{equation}
  \label{eq:uc:objective}%
  \min_{ \substack{ \mathbf{x},\mathbf{y},\mathbf{z}, \\
  \mathbf{P}_{t,g},\mathbf{V}_{t},\mathbf{q}}} \qquad \sum_{t=1}^{T} \sum_{g=1}^{N_{G}}
  \mathbf{C}_{g} \bullet \mathbf{P}_{t,g} +
    \mathbf{c}^{\intercal} \mathbf{x} +\mathbf{u}^{\intercal} \mathbf{y} +
    \mathbf{h}^{\intercal} \mathbf{z}
\end{equation}
subject to:
\begin{eqnarray}
  \label{eq:uc:active-power-balance}%
  \sum_{\forall g\in\mathcal{G}_{i}} \left( \mathbf{A} \bullet \mathbf{P}_{t,g}+ \underline{p}_{g} \cdot x_{t,g} \right)+\mathbf{Y}_{i}\bullet \mathbf{V}_{t}&=&dp_{t,i}\\
  \label{eq:uc:reactive-power-balance}%
  \sum_{\forall g\in\mathcal{G}_{i}} \left(\overset{\Delta}{q}_{t,g}+ \underline{q}_{g} \cdot x_{t,g} \right)+\widetilde{\mathbf{Y}}_{i}\bullet \mathbf{V}_{t}&=&dq_{t,i}\\
  \label{eq:uc:active-power}%
  \mathbf{A} \bullet \mathbf{P}_{t,g} - \left( \overline{p}_{g} - \underline{p}_{g}\right) x_{t,g} &\leqslant& 0\\
    \label{eq:uc:reactive-power}%
  \overset{\Delta}{q}_{t,g} - \left( \overline{q}_{g} - \underline{q}_{g}\right) x_{t,g} &\leqslant& 0\\
  \label{eq:uc:startup}%
  x_{t,g} - \textbf{1}_{(t>1)} \cdot x_{t-1,g} - \textbf{1}_{(t=1)} \cdot x_{0,g}&\leqslant& y_{t,g} \\
  \label{eq:uc:shutdown}%
  \textbf{1}_{(t>1)} \cdot x_{t-1,g} + \textbf{1}_{(t=1)} \cdot x_{0,g} - x_{t,g} &\leqslant& z_{t,g}\\
  \label{eq:uc:time-up}%
  \omega_{(t,T_{on})}\left(x_{t}- \textbf{1}_{(t>1)}\cdot x_{t-1}\right) - \sum_{j=1}^{\omega_{(t,T_{on})}} x_{j+t-1,g} &\leqslant& 0\\
  \label{eq:uc:time-down}%
  \omega_{(t,T_{off})}\left(\textbf{1}_{(t>1)}\cdot x_{t-1} - x_{t}\right) - \sum_{j=1}^{\omega_{(t,T_{off})}}\left(1-x_{j+t-1,g}\right) &\leqslant& 0\\
  \label{eq:uc:spinning-reserve}%
  \sum_{g=1}^{N_{G}} \left( \mathbf{A} \bullet \mathbf{P}_{t,g} - \left( \overline{p}_{g}-\underline{p}_{g} \right) \cdot x_{t,g} \right) &\leqslant& - SR_{t}\\
  \label{eq:uc:ramp-up}%
    \left(\mathbf{A}\bullet\mathbf{P}_{t,g}+\underline{p}_{g}\cdot x_{t,g} \right) - \textbf{1}_{(t=1)}\cdot  p_{0,g}
    -\textbf{1}_{(t>1)} \left(\mathbf{A}\bullet\mathbf{P}_{t-1,g} + \underline{p}_{g}\cdot x_{t-1,g}\right) &\leqslant& RU_{g}\\
      \label{eq:uc:ramp-down}%
  \textbf{1}_{(t=1)}\cdot  p_{0,g} - \left( \mathbf{A} \bullet \mathbf{P}_{t,g} + \underline{p}_{g} \cdot x_{t,g} \right)
    + \textbf{1}_{(t>1)} \left( \mathbf{A} \bullet \mathbf{P}_{t-1,g} + \underline{p}_{g}\cdot x_{t-1,g}\right) &\leqslant& RD_{g}\\
  \label{eq:uc:flow-limits}%
  -\overline{F}_{i,j} \;\;\, \leqslant \;\;\, \mathbf{Y}_{i,j} \bullet \mathbf{V}_{t} &\leqslant& \overline{F}_{i,j}\\
  \label{eq:uc:voltage}%
  \underline{V}_{i}^{2} \;\;\, \leqslant \;\;\, \mathbf{E}_{i}\bullet \mathbf{V}_{t} &\leqslant& \overline{V}_{i}^{2}\\
  \label{eq:uc:semidefinite}%
  \mathbf{P}_{t,g}, \, \mathbf{V}_{t} &\succeq& 0\\
  \label{eq:uc:positive}%
  \mathbf{q} &\geq& 0\\
  \label{eq:uc:rank}%
  \textrm{rank}(\mathbf{V}_{t}) &=& 1\\
  \label{eq:uc:binary}%
  \mathbf{x},\mathbf{y},\mathbf{z}  &\in& \left\{ 0,1 \right \}^{T \cdot N_{G}} \quad
\end{eqnarray}

The UC problem~\eqref{eq:uc:objective}--\eqref{eq:uc:binary} is defined over a variable space consisting of active and reactive power, bus voltages, and discrete decisions describing unit status, and startup/shutdown events. Besides the constraints associated with AC power flow, the constraint set includes minimum up and down times, ramp, and spinning reserve inequalities. Its formulation enables the application of a Benders' decomposition technique~\cite{Benders} by means of a convex relaxation of the rectangular AC OPF problem based on~\cite{5971792}, that seeks to provide optimal solutions to active power generation and UC decisions subject to AC power flow constraints over the voltage variable space.

\subsection{Problem dimensions}
\label{sec:problem:dimensions}
The UC problem is defined over a set of $g=1,2,\ldots,N_{G}$ generators, $i=1,2,\ldots,N_{B}$ buses, and discrete time stages $t=1,2,\ldots,T$, most commonly in hourly steps, and lines $l=1,2,\ldots,N_{L}$ for every undirected pair $(i,j)$ of buses $1\leqslant i \leqslant N_{B}$, and $1 \leqslant j \leqslant N_{B}$, $i\neq j$ for which there is a connecting transmission line.

\subsection{Problem variables}%
\label{sec:problem:variables}

\subsubsection{Startup/shutdown}
\label{sec:problem:variables:binary}
Binary vectors $\mathbf{x}$, $\mathbf{y}$, and $\mathbf{z}$ represent up/down statuses, and startup and shutdown events, respectively, for a given generator $g$ at every $t$. If generator $g$ is up-and-running during $t$ then $x_{t,g}$ equals $1$, otherwise it is $0$. On the other hand, if generator $g$ is brought up at time $t$, then $y_{t,g}$ equals $1$, or $0$ otherwise. Analogously, if generator $g$ is shutdown at time $t$, then $z_{t,g}$ equals $1$, or $0$
otherwise. The indicator function $\textbf{1}_{(t>1)}$ is
equal to $1$ for time steps $t>1$, or $0$ otherwise. The value of $x_{0,g}$ is assumed to be known. Therefore:
\begin{eqnarray*}
  \mathbf{x} &\doteq& \left[ x_{1,1},\ldots,x_{1,N_{G}},\ldots,x_{T,1},\ldots,x_{T,N_{G}}\right]^{\intercal},
  \mathbf{y} \doteq \left[ y_{1,1},\ldots,y_{1,N_{G}},\ldots,y_{T,1},\ldots,y_{T,N_{G}}\right]^{\intercal}, \textrm{ and}\\
  \mathbf{z} &\doteq& \left[ z_{1,1},\ldots,z_{1,N_{G}},\ldots,z_{T,1},\ldots,z_{T,N_{G}}\right]^{\intercal}.
\end{eqnarray*}

\subsubsection{Power generation}
\label{sec:problem:variables:power}
Active power generation $p_{t,g}$ is subject to minimum $\underline{p}_{g}$ and maximum $\overline{p}_{g}$ limits. Let:
\begin{displaymath}
  \overset{\Delta}{p}_{t,g} \doteq p_{t,g} - \underline{p}_{g}
\end{displaymath}
define the incremental active power generation if generator $g$ is up-and-running. Then $\overset{\Delta}{p}_{t,g}$ can be alternatively written as a Frobenius product $\mathbf{A} \bullet \mathbf{P}_{t,g}$, where:
\begin{equation*}
  \mathbf{A} \doteq \left[ \begin{array}{cc}
    0 & 1/2 \\
    1/2 & 0
  \end{array} \right], \textrm{ and} \;
  \mathbf{P}_{t,g}  \doteq \left[ \begin{array}{cc}
      \overset{\Delta}{p}_{t,g}^2 & \overset{\Delta}{p}_{t,g} \\
      \overset{\Delta}{p}_{t,g} & 1
    \end{array} \right],
\end{equation*}
and the Frobenius product of two matrices is defined as the sum of element-wise products, such that:
\begin{displaymath}
  \overset{\Delta}{p}_{t,g} = \mathbf{A} \bullet \mathbf{P}_{t,g} = \frac{1}{2} \overset{\Delta}{p}_{t,g} + \frac{1}{2} \overset{\Delta}{p}_{t,g}.
\end{displaymath}
Similarly, reactive power generation is represented as the difference with the minimum reactive power so:
\begin{displaymath}
  \overset{\Delta}{q}_{t,g} \doteq q_{t,g} - \underline{q}_{g}
\end{displaymath}
\begin{equation*}
  \mathbf{q} \doteq \left[ \overset{\Delta}{q}_{1,1},\ldots,\overset{\Delta}{q}_{1,N_{G}},\ldots,\overset{\Delta}{q}_{T,1},\ldots,\overset{\Delta}{q}_{T,N_{G}}\right]^{\intercal}.
\end{equation*}

\subsubsection{Voltage}
\label{sec:problem:variables:voltage}
Real $e_{t,i}$ and imaginary $f_{t,i}$ parts of voltage at bus $1\leqslant i\leqslant N_{B}$ define the voltage vector as follows:
\begin{displaymath}
  \mathbf{v}_{t} \doteq \left[ e_{t,1},\ldots,e_{t,N_{B}},f_{t,1},\ldots,f_{t,N_{B}}  \right]^{\intercal} \; \in \; \mathbb{R}^{2\cdot N_{B}},
\end{displaymath}
such that the voltage variable $\mathbf{V}_{t} \in \mathbb{R}^{2\cdot N_{B} \times 2\cdot N_{B}}$ is given by the outer product of the voltage vector:
\begin{equation}
  \label{eq:problem:variables:voltage:vv}
  \mathbf{V}_{t} = \mathbf{v}_{t}\mathbf{v}_{t}^{\intercal}.
\end{equation}

\subsection{Objective function}
\label{sec:problem:objective}
The objective function in~\eqref{eq:uc:objective} represents costs associated with active power generation, as well as generator startup and shutdown events.
Let the generation costs of an up-and-running generator $g$ during step $t$ be approximated by a quadratic function of the form:
\begin{displaymath}
  \alpha_{g} \cdot p_{t,g}^{2} + \beta_{g} \cdot p_{t,g} + \gamma_{g},
\end{displaymath}
then it can alternatively be expressed in terms of incremental power generation:
\begin{displaymath}
  \alpha_{g} \cdot \left(\overset{\Delta}{p}_{t,g} + \underline{p}_{t,g}\right)^2 + \beta_{g} \cdot \left(\overset{\Delta}{p}_{t,g} + \underline{p}_{t,g}\right) + \gamma_{g}.
\end{displaymath}

In~\eqref{eq:uc:objective} the incremental costs are expressed as a Frobenius product $\mathbf{C}_{g} \bullet \mathbf{P}_{t,g}$, such that:
\begin{displaymath}
  \mathbf{C}_{g} \doteq \left[ \begin{array}{cc}
      \alpha_{g} & \alpha_{g}\cdot \underline{p}_{g} + \beta_{g}/2\\
      \alpha_{g}\cdot \underline{p}_{g} + \beta_{g}/2 & 0
    \end{array} \right]
\end{displaymath}
and the minimum generation costs are given by the inner product $\mathbf{c}^{\intercal} \mathbf{x}$, such that:
\begin{displaymath}
  \mathbf{c} \doteq \left[c_{1,1}, \ldots ,c_{1,N_{G}},\ldots , c_{T,1}, \ldots ,c_{T,N_{G}} \right]^{\intercal} \; \in \; \mathbb{R}_{+}^{T \cdot N_{G}}
\end{displaymath}
where:
\begin{displaymath}
    c_{t,g} \doteq \alpha_{g} \cdot \underline{p}_{g}^{2}+\beta_{g}\cdot \underline{p}_{g}+\gamma_{g},
\end{displaymath}
is constant over time, and $\alpha_{g}$, $\beta_{g}$, and $\gamma_{g}$ are the coefficients of the quadratic active power generation cost function for generator $g$.

Finally, vectors $\mathbf{u}, \mathbf{h} \in \mathbb{R}_{+}^{T \cdot N_{G}}$ represent cost coefficients of startup and shutdown events for every generator $g$ and time $t$, respectively.

\subsection{Power balance constraints}
\subsubsection{Active power balance}
\label{sec:problem:constraints:active-power-balance}
In Eq.~\eqref{eq:uc:active-power-balance} active power balance constraints are expressed as the difference between active load demand $dp_{t,i}$ and the sum of total active power generation and net active power flow to bus $i=1,2,\ldots,N_{B}$ in every $t=1,2,\ldots,T$. Active power generation is expressed as the sum of incremental $\mathbf{A} \bullet \mathbf{P}_{t,g}$ and minimum generation $\underline{p}_{g} \cdot x_{t,g}$ at all generators $g$ connected to bus $i$, as represented by set $\mathcal{G}_{i}$. Net active power flow injection to bus $i$ is given by $\mathbf{Y}_{i}\bullet \mathbf{V}_{t}$, such that:
\begin{displaymath}
  \mathbf{Y}_{i} \doteq -\sum_{\forall j \in \Omega_{i}} \mathbf{Y}_{i,j}
\end{displaymath}
where $\Omega_{i}$ is the set of bus indexes of all buses directly connected to bus $i$, so that the pair $(i,j)$ exists, and:
\begin{displaymath}
  \mathbf{Y}_{i,j} \doteq \frac{1}{2} \left[ \begin{array}{cc}
    G_{i,j}+G_{i,j}^{\intercal} & B_{i,j}^{\intercal}-B_{i,j} \\
    B_{i,j}-B_{i,j}^{\intercal} & G_{i,j}+G_{i,j}^{\intercal}
  \end{array} \right] \; \in \; \mathbb{R}^{ (2\cdot N_{B}) \times (2\cdot N_{B})}
\end{displaymath}
is the admittance matrix, and:
\begin{equation*}
  G_{i,j} \doteq g_{i,j} \left( \boldsymbol{\xi}_{i}\boldsymbol{\xi}_{i}^{\intercal}- \boldsymbol{\xi}_{i}\boldsymbol{\xi}_{j}^{\intercal} \right), \;
  B_{i,j} \doteq b_{i,j} \left( \boldsymbol{\xi}_{i}\boldsymbol{\xi}_{i}^{\intercal}- \boldsymbol{\xi}_{i}\boldsymbol{\xi}_{j}^{\intercal}  \right),
\end{equation*}
where $g_{i,j}$ and $b_{i,j}$ represent line conductance and susceptance, respectively, and $\boldsymbol{\xi}_{i}\in\mathbb{R}^{N_{B}}$ is the standard basis vector.

\subsubsection{Reactive power balance}
\label{sec:problem:constraints:reactive-power-balance}
In Eq.~\eqref{eq:uc:reactive-power-balance} reactive power balance constraints are expressed as the difference between reactive load demand $dq_{t,i}$ and the sum of total active power generation and net reactive power flow to bus $i=1,2,\ldots,N_{B}$ in every $t=1,2,\ldots,T$. Reactive power generation is expressed as the sum of incremental $\overset{\Delta}{q}_{t,g} $ and minimum generation $\underline{q}_{g} \cdot x_{t,g}$ at all generators $g$ connected to bus $i$. Net reactive power flow injection to bus $i$ is given by $\widetilde{\mathbf{Y}}_{i}\bullet \mathbf{V}_{t}$, such that:
\begin{displaymath}
  \widetilde{\mathbf{Y}}_{i} \doteq -\sum_{\forall j \in \Omega_{i}} \widetilde{\mathbf{Y}}_{i,j}
\end{displaymath}
and:
\begin{displaymath}
 \widetilde{\mathbf{Y}}_{i,j} \doteq \frac{1}{2} \left[ \begin{array}{cc}
    B_{i,j}+B_{i,j}^{\intercal} & G_{i,j}-G^{\intercal}_{i,j} \\
    G_{i,j}^{\intercal}-G_{i,j} & B_{i,j}+B_{i,j}^{\intercal}
  \end{array} \right] \; \in \; \mathbb{R}^{ (2\cdot N_{B}) \times (2\cdot N_{B})}
\end{displaymath}

\subsection{Active power constraints}
\label{sec:problem:constraints:active-reactive-power}
Maximum active and reactive power generation constraints are expressed in Eq.~\eqref{eq:uc:active-power} and Eq.~\eqref{eq:uc:reactive-power} respectively, for $g=1,2,\ldots,N_{G}$, and $t=1,2,\ldots,T$ in terms of incremental power generation.

\subsection{Startup and shutdown constraints}
\label{sec:problem:constraints:startup-shutdown}
In Eq.~\eqref{eq:uc:startup} startup constraints are represented for every $g=1,2,\ldots,N_{G}$, and $t=1,2,\ldots,T$. Analogously, shutdown constraints are represented in Eq.~\eqref{eq:uc:shutdown} along the same dimensions.

\subsection{Minimum up and down time constraints}
\label{sec:problem:constraints:minimum-up-down}
In Eq.~\eqref{eq:uc:time-up}, and Eq.~\eqref{eq:uc:time-down} minimum up $T_{on}$ and down $T_{off}$ time constraints are represented, respectively, for every $g=1,2,\ldots,N_{G}$, and $t=1,2,\ldots,T$, where the sum of time steps (second term) in up (down) time after a startup (shutdown) event occurs (first term) is equal to the minimum number of time steps $\omega_{(t,\tau)}$ required for a unit to be up (down). This formulation is an adaptation of that presented in~\cite{4784461}. Constant value $T_{0}$ represents a unit status immediately before the beginning of the considered time horizon---it is positive if unit was up, or negative otherwise, where:
\begin{displaymath}
  \omega_{(t,\tau)} = \left\{ \begin{array}{lr}
    \min \left(\tau-T_{0}, \tau \right) & : t=1,\\
    \min \left(\tau, T-t+1 \right) & : t>1.
  \end{array} \right.
\end{displaymath}

\subsection{Spinning reserve constraints}
\label{sec:problem:constraints:reserve}
In Eq.~\eqref{eq:uc:spinning-reserve} constraints on residual active power generation availability $SR_{t}$, also known as spinning reserve, are imposed for every $t=1,2,\ldots,T$.

\subsection{Maximum ramp up and down rate constraints}
\label{sec:problem:constraints:ramp}
Active power generation is subject to maximum ramp up $RU_{g}$ and down $RD_{g}$ rates between consecutive time steps, as expressed by Eq.~\eqref{eq:uc:ramp-up}, and Eq.~\eqref{eq:uc:ramp-down}, respectively, for every $g=1,2,\ldots,N_{G}$, and $t=1,2,\ldots,T$. Where $p_{0,g}$ represents active power dispatch and its current value is assumed to be known. 

\subsection{Power flow constraints}
\label{sec:problem:constraints:power-flow}
In Eq.~\eqref{eq:uc:flow-limits} maximum power flow constraints $\overline{F}_{i,j}$ are imposed on power flow $\mathbf{Y}_{i,j} \bullet \mathbf{V}_{t}$ in both directions of every line $l=1,2,\ldots,N_{L}$ such that $(i,j)$ exists.

\begin{remark}
	This constraint represents a limit on active power flow in the network, and does not cover limits on reactive power flow.
\end{remark}

\subsection{Voltage constraints}
\label{sec:problem:constraints:voltage}
The module of the voltage is subject to minimum $\underline{V}_{i}$ and maximum $\overline{V}_{i}$ limits in every $i=1,2,\ldots,N_{B}$, as expressed in Eq.~\eqref{eq:uc:voltage}, where:
\begin{displaymath}
  \mathbf{E}_{i} \doteq \left[ \begin{array}{cl}
    \boldsymbol{\xi}_{i}\boldsymbol{\xi}_{i}^{\intercal} & \mathbf{0}\\
    \mathbf{0} & \boldsymbol{\xi}_{i}\boldsymbol{\xi}_{i}^{\intercal}
  \end{array} \right].
\end{displaymath}

\subsection{Rank constraints}
\label{sec:problem:constraints:rank}
In Eq.~\eqref{eq:uc:rank} a rank-1 constraint is imposed on $\mathbf{V}_{t}$ as a consequence of the transformation in Eq.~\eqref{eq:problem:variables:voltage:vv}, and the relaxation $\mathbf{V}_{t} - \mathbf{v}\mathbf{v}^{\intercal} \succeq 0$ implicit to Eq.~\eqref{eq:uc:semidefinite}.

On the other hand, if matrix $\mathbf{P}_{t,g} \succeq 0$, and therefore, its terms $P_{t,g,1,2}=\overset{\Delta}{p}_{t,g}$ and $P_{t,g,2,2}$ equal $1$, then the following holds for the determinant of a semidefinite matrix:
\begin{equation*}
  P_{t,g,1,1} \cdot 1 - P_{t,g,1,2} \cdot P_{t,g,2,1}  \geqslant 0 , \;
  P_{t,g,1,1}  \geqslant \overset{\Delta}{p}_{t,g}^2
\end{equation*}

\begin{remark}
If $\alpha>0$, which can be reasonably assumed since thermoelectric power generation can be satisfactorily approximated by convex functions, we can deduce that the relationship $\mathbf{C}_{g} \bullet \mathbf{P}_{t,g} \geqslant \alpha_{t,g} \cdot p_{t,g}^2 + \beta_{t,g} \cdot p_{t,g}$ is implicit to $\mathbf{P}_{t,g} \succeq 0$, and, in a minimization algorithm, this relationship will be active since $P_{t,g,1,1}$ does not have any other lower bound. Furthermore, from such observation we can also deduce that in the optimal solution matrix $\mathbf{P}_{t,g}$ is rank-$1$.
\end{remark}

\section{Solution Methodology}%
\label{sec:solution}
In order to solve the mixed-integer SDP relaxation of the UC problem defined in~\eqref{eq:uc:objective}--\eqref{eq:uc:binary}, we propose the use of a generalized Benders'-like decomposition scheme. Initially, a feasible UC schedule is obtained with a relaxed bus voltage active power dispatch. In~\cite{5971792}, the OPF problem is numerically perturbed in order to yield rank-2 solutions, whereas in this paper we choose to use a rank reduction procedure (RRP) that avoids such problem-dependent perturbation, and also provides low-rank solutions to the voltage variables with optimal active power generation and UC dispatches.

In the algorithm presented in the following sections, a decomposition scheme for polynomial optimization problems~\cite{benders-SDP} is adapted according to an application of the generalized Benders' technique for mixed-integer variables~\cite{Benders} to the UC problem, therefore rendering a Benders'-like decomposition method for mixed-integer semidefinite programming (MI-SDP) problems. The method is summarized in Fig.~\ref{fig:scheme}.

\tikzstyle{decision} = [diamond, aspect=2, draw, text width=4.5em, text badly centered, node distance=2cm, inner sep=0pt]
\tikzstyle{block} = [rectangle, draw, text centered, rounded corners, minimum height=4em]
\tikzstyle{line} = [draw, -latex']

\begin{figure}[!h]
  \begin{center}
  \begin{scriptsize}
  \begin{tikzpicture}[node distance = 1.8cm, auto]
    \node [block] (init) {Solve Master Problem $MM^{(k)}$};
    \node [decision, below of=init] (while) {Did it converge?};
    \node [block,text width=7em, right of=while, node distance = 3cm] (fim) {Infeasible UC};
    \node [block, below of=while] (minim) {Solve Subproblem $S^{(k)}$};
    \node [decision, below of=minim] (desition) {Is it feasible?};
    \node [block,text width=7em, left of=desition, node distance = 3cm ] (feasible) {Add feasibility cut~\eqref{eq:feasibility}};
    \node [decision, text width=8em,below of=desition] (mismatch) {$UB^{(k)}=LB^{(k)}$};
    \node [block,text width=7em, left of=mismatch, node distance = 3cm] (optimality) {Add optimality cut~\eqref{eq:optimality}};
        \node [block,text width=8em, below of=mismatch] (rank) { Matrix $\mathbf{V}^{*}_{k}$ Rank Reduction (Section~\ref{sec:rank})};
        \node [block,text width=5em, below of=rank] (end) { End };
    \path [line] (init) -- (while);
    \path [line] (while) --node {Yes} (minim);
    \path [line] (while) --node {No} (fim);
    \path [line] (minim) -- (desition);
    \path [line] (desition) -- node {Yes} (mismatch);
    \path [line] (desition) -- node {No} (feasible);
    \path [line] (mismatch) -- node {Yes} (rank);
    \path [line] (rank) -- (end);
    \path [line] (mismatch) -- node {No} (optimality);
	\path [line] (feasible.west) to [out=120, in=180] (init.west);
        \path [line] (optimality.west) to [out=120, in=180] (init.west);
      \end{tikzpicture}
    \end{scriptsize}
  \end{center}
  \caption{Flowchart summarizing the procedures of the solution methodology proposed for solution of the UC problem. It includes the Benders'-like decomposition scheme with solution of the mixed-integer master problem and SDP relaxation of the optimal power flow subproblem, inclusion of feasibility and optimality cuts, and voltage matrix rank reduction procedure.}
  \label{fig:scheme}
\end{figure}
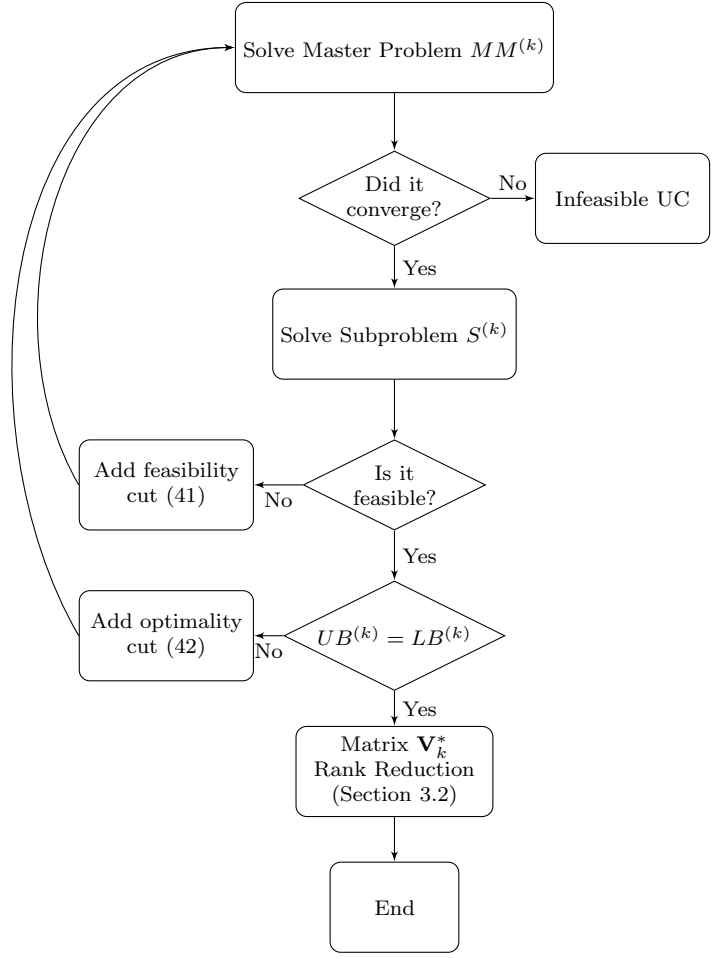

\subsection{Benders'-like decomposition}
\label{sec:solution:benders}
We can write the UC problem~\eqref{eq:uc:objective}--\eqref{eq:uc:binary} in matrix notation in order to simplify the development of our Benders' decomposition scheme. With that goal in mind, let $\mathbf{C}\in\mathbb{R}^{(2\cdot T \cdot N_{G})\times (2 \cdot T \cdot N_{G})}$ represent the incremental cost matrix for all generators and time steps, such that:
\begin{displaymath}
  \mathbf{C} \doteq \textrm{diag}\left(\mathbf{C}_{1},\ldots,\mathbf{C}_{N_{G}},\ldots, \mathbf{C}_{1},\ldots,\mathbf{C}_{N_{G}}\right).
\end{displaymath}
Consequently, we need to represent incremental active power generation variables $\mathbf{P} \in \mathbb{R}^{(2 \cdot T \cdot N_{G}) \times( 2 \cdot T \cdot N_{G})}$ as a block-diagonal matrix along the same dimensions, with $\mathbf{P}_{t,g} \in \mathbb{R}^{2 \times 2}$ in its diagonal, such that:
\begin{displaymath}
  \mathbf{P} \doteq \textrm{diag}\left(\mathbf{P}_{1,1},\ldots,\mathbf{P}_{1,g}, \ldots, \mathbf{P}_{T,1},\ldots,\mathbf{P}_{T,N_{G}}\right).
\end{displaymath}
A voltage matrix $\mathbf{V} \in \mathbb{R}^{(2\cdot T \cdot N_{B}) \times (2 \cdot T \cdot N_{B})}$ is constructed in analogous terms:
\begin{displaymath}
  \mathbf{V} \doteq \textrm{diag}\left(\mathbf{V}_{1}, \mathbf{V}_{2}, \ldots, \mathbf{V}_{T}\right).
\end{displaymath}

Then the resulting equivalent UC problem in matrix notation with relaxation of the rank constraints is expressed as follows:
\begin{equation}
  \label{eq:uc:m:objective}%
  \displaystyle \min_{ \substack{ \mathbf{x},\mathbf{y},\mathbf{z}, \\
  \mathbf{P},\mathbf{V},\mathbf{q}}} \qquad
  \mathbf{C} \bullet \mathbf{P} +
    \mathbf{c}^{\intercal} \mathbf{x} +\mathbf{u}^{\intercal} \mathbf{y} +
    \mathbf{h}^{\intercal} \mathbf{z} \qquad\qquad
\end{equation}
\begin{eqnarray}
  \label{eq:uc:m:active-power-balance}%
  \text{s. t. }\;
  \mathcal{A}\mathbf{1}(\mathbf{P}) + \underline{\mathbf{P1}} \mathbf{x} + \mathcal{Y}\mathbf{1}(\mathbf{V})
    &=& \mathbf{d}p\\
   \label{eq:uc:m:reactive-power-balance}%
  \mathbf{q} + \underline{\mathbf{Q}} \mathbf{x} + \mathcal{Y}\mathbf{2}(\mathbf{V})
    &=& \mathbf{d}q\\
  \label{eq:uc:m:active-power}%
  \mathcal{A} \mathbf{2}(\mathbf{P}) - \overline{\mathbf{P1}} \mathbf{x}
    &\leqslant& \mathbf{0}\\
  \label{eq:uc:m:reactive-power}%
  \mathbf{q} -\overline{\mathbf{Q}}\mathbf{x}  &\leqslant& \mathbf{0}\\
  \label{eq:uc:m:startup}%
  \mathbf{Ux} - \mathbf{y}
    &\leqslant& \mathbf{a1}\\
  \label{eq:uc:m:shutdown}%
  -\mathbf{Ux} - \mathbf{z}
    &\leqslant& - \mathbf{a1}\\
  \label{eq:uc:m:time-up}%
  \mathbf{M1x}
    &\leqslant& \mathbf{0}\\
  \label{eq:uc:m:time-down}%
  \mathbf{M2x}
    &\leqslant& \mathbf{a2}
\end{eqnarray}

\begin{eqnarray}
   \label{eq:uc:m:spinning-reserve}%
  \mathcal{A}\mathbf{3}(\mathbf{P}) + \overline{\mathbf{P2}}\mathbf{x}
    &\leqslant& \mathbf{r}\\
  \label{eq:uc:m:ramp-up}%
  \mathcal{A}\mathbf{4}(\mathbf{P})+\overline{\mathbf{P3}}\mathbf{x}
    &\leqslant& \mathbf{g1}\\
  \label{eq:uc:m:ramp-down}%
  -\mathcal{A} \mathbf{4}(\mathbf{P})-\overline{\mathbf{P3}}\mathbf{x}
    &\leqslant& \mathbf{g2}\\
  \label{eq:uc:m:flow-limits}%
  \mathcal{Y}\mathbf{3}(\mathbf{V})
    &\leqslant& \mathbf{f}\\
  \label{eq:uc:m:voltage}%
  \mathcal{Y} \mathbf{4}(\mathbf{V})
    &\leqslant& \mathbf{b}\\
  \label{eq:uc:m:semidefinite}%
  \mathbf{P}, \, \mathbf{V} &\succeq& 0\\
  \label{eq:uc:m:positive}%
  \mathbf{q} &\geq& 0\\
  \label{eq:uc:m:rank}%
  \textrm{rank}(\mathbf{V}) &=& 1\\
  \label{eq:uc:m:binary}%
  \mathbf{x},\mathbf{y},\mathbf{z}  &\in& \left\{ 0,1 \right \}^{T \cdot N_{G}}
\end{eqnarray}

Let $\boldsymbol{\varrho}_{(t,g)}\doteq\boldsymbol{\xi}_{N_{G} \cdot \left(t-1 \right)+g}$ define a standard basis vector of appropriate dimension.

Active power balance constraints~\eqref{eq:uc:active-power-balance} are expressed in Eq.~\eqref{eq:uc:m:active-power-balance} by linear mappings $\mathcal{A} \mathbf{1} (\cdot)$ and $\mathcal{Y}\mathbf{1}\left( \cdot \right)$, representing incremental active power generation and net power flow injection, respectively. They are composed of matrices $\mathbf{A1}_{r} \in \mathbb{R}^{(2 \cdot T \cdot N_{G}) \times (2 \cdot T \cdot N_{G} )}$
and $\mathbf{Y1}_{r} \in \mathbb{R}^{(2\cdot N_{B} )\times ( 2\cdot N_{B})}$, respectively:
\begin{eqnarray*}
  \displaystyle\mathbf{A1}_{r} &\doteq& \sum_{g \in \mathcal{G}_{i}} \boldsymbol{\varrho}_{(t,g)}\boldsymbol{\varrho}_{(t,g)}^{\intercal} \otimes \mathbf{A},\\
  \mathbf{Y1}_{r} &\doteq& \mathbf{Y}_{i}.
\end{eqnarray*}
where $r=N_{B}\cdot(t-1)+i$, and $\otimes$ is the Kronecker product. In addition, matrix $\underline{ \mathbf{P1}} \in \mathbb{R}^{ (T \cdot N_{B}) \times( T \cdot N_{G} ) }$ represents minimum active power generation, such that:
\begin{displaymath}
  \underline{ \mathbf{ P1}} \doteq \sum_{t=1}^{T} \sum_{i=1}^{N_{B}} \sum_{\forall g \in \mathcal{G}_{i}} \boldsymbol{\xi}_{r} \boldsymbol{\varrho}_{(t,g) }^{\intercal},
\end{displaymath}
and, finally:
\begin{displaymath}
  \mathbf{d}p \doteq \left[ d_{1,1}, \ldots, d_{1,N_{B}}, \ldots, d_{T,1}, \ldots, d_{T,N_{B}} \right]^{\intercal} \in \Real^{T\cdot N_{B}}.
\end{displaymath}

Reactive power balance constraints~\eqref{eq:uc:reactive-power-balance} are expressed in Eq.~\eqref{eq:uc:m:reactive-power-balance} by vectors $\mathbf{q}$, $\underline{\mathbf{q}}$  representing incremental and minimum reactive power generation respectively. Vector $\underline{\mathbf{q}}$ is composed by:
\begin{equation}
\underline{\mathbf{Q}} = \text{diag} \left(\underline{q}_{1},\ldots,\underline{q}_{N_{G}},\underline{q}_{1},\ldots,\underline{q}_{N_{G}} \right)\in \mathbb{R}^{T\cdot N_{G}}
\end{equation}
and linear mapping $\mathcal{Y}\mathbf{2}(\cdot)$ representing minimum reactive power generation and net power flow injection, is composed by matrices $\mathbf{Y2}_{r}$:
\begin{equation}
\mathbf{Y2}_{r} \doteq \widetilde{\mathbf{Y}}_{i}.
\end{equation}
where $r=N_{B}\cdot(t-1)+i$, also the reactive load demand is:
\begin{displaymath}
  \mathbf{d}q \doteq \left[ dq_{1,1}, \ldots, dq_{1,N_{B}}, \ldots, dq_{T,1}, \ldots, dq_{T,N_{B}} \right]^{\intercal} \in \Real^{T\cdot N_{B}}.
\end{displaymath}

Active power constraints~\eqref{eq:uc:active-power} are expressed in Eq.~\eqref{eq:uc:m:active-power} in terms of linear mapping $\mathcal{A} \mathbf{2} (\cdot)$, which is composed of matrices $\mathbf{A2}_{r} \in \mathbb{R}^{(2 \cdot T \cdot N_{G}) \times (2 \cdot T \cdot N_{G}) }$,
and $\overline{ \mathbf{ P1}} \in \mathbb{R}^{ (T \cdot N_{B}) \times( T \cdot N_{G} ) }$:
\begin{eqnarray*}
  \mathbf{A2}_{r} &\doteq& \boldsymbol{\varrho}_{(t,g)} \boldsymbol{\varrho}_{(t,g)}^{\intercal} \otimes \mathbf{A}\\
  \displaystyle\overline{ \mathbf{ P1}} &\doteq& \sum_{t=1}^{T} \overline{\Delta P}_{t,g} \cdot \boldsymbol{\varrho}_{(t,g)} \boldsymbol{\varrho}_{(t,g)}^{\intercal}
\end{eqnarray*}

Reactive power constraints~\eqref{eq:uc:reactive-power} are expressed in Eq.~\eqref{eq:uc:m:reactive-power} with upper bound vector:
\begin{equation}
\overline{\mathbf{Q}} = \text{diag} \left( \overline{q}_{1}-\underline{q}_{1},\ldots, \overline{q}_{N_{G}}-\underline{q}_{N_{G}},\overline{q}_{1}-\underline{q}_{1},\ldots, \overline{q}_{N_{G}}-\underline{q}_{N_{G}}\right)
\end{equation}

Startup~\eqref{eq:uc:startup} and shutdown~\eqref{eq:uc:shutdown} contraints are expressed in Eqs.~\eqref{eq:uc:m:startup} and~\eqref{eq:uc:m:shutdown}, respectively, in terms of $\mathbf{U} \in \mathbb{R}^{(T \cdot N_{G}) \times (T \cdot N_{G})}$ and $\mathbf{a1} \in \mathbb{R}^{T \cdot N_{G}}$:
\begin{eqnarray*}
  \mathbf{U} &\doteq& \mathbf{I}- \sum_{t=2}^{T} \sum_{g=1}^{N_{G}} \boldsymbol{\varrho}_{(t,g)}\boldsymbol{\varrho}_{(t-1,g)}^{\intercal}\\
  \mathbf{a1} &\doteq& \sum_{g=1}^{N_{G}} x_{0,g} \cdot \boldsymbol{\varrho}_{(t,g)}
\end{eqnarray*}

Minimum up~\eqref{eq:uc:time-up} and down~\eqref{eq:uc:time-down} time constraints are expressed in Eqs.~\eqref{eq:uc:m:time-up} and~\eqref{eq:uc:m:time-down}, respectively, in terms of matrices $\mathbf{M1},  \mathbf{M2} \in \mathbb{R}^{(T\cdot N_{G}) \times (T\cdot N_{G})}$,
such that:
\begin{displaymath}
    \sum_{t=1}^{T}\sum_{g=1}^{N_{G}}  \left(  \omega _{\left(t, T_{on} \right)} \left( \boldsymbol{\varrho}_{(t,g)} \boldsymbol{\varrho}_{(t,g)}^{\intercal} - \textbf{1}_{(t>1)} \cdot \boldsymbol{\varrho}_{(t-1,g)} \boldsymbol{\varrho}_{(t-1,g)}^{\intercal} \right)  \right.
    \left. - \sum_{j=1}^{\omega_{(t,T_{on})} }   \boldsymbol{\varrho}_{(j+t-1,g)} \boldsymbol{\varrho}_{(j+t-1,g)}^{\intercal}   \right) = \mathbf{M1}
\end{displaymath}
and:
\begin{displaymath}
    \sum_{t=1}^{T}\sum_{g=1}^{N_{G}} \left(\omega _{\left(t, T_{off} \right) } \left(  \textbf{1}_{(t>1)} \cdot \boldsymbol{\varrho}_{(t-1,g)} \boldsymbol{\varrho}_{(t,g)}^{\intercal} -  \boldsymbol{\varrho}_{(t-1,g)}  \boldsymbol{\varrho}_{(t-1,g)}^{\intercal} \right) \right.
    \left. + \sum_{j=1}^{\omega_{(t,T_{off})} }   \boldsymbol{\varrho}_{(j+t-1,g)} \boldsymbol{\varrho}_{(j+t-1,g)}^{\intercal} \right) = \mathbf{M2}
\end{displaymath}
where $\mathbf{a2} \in \mathbb{R}^{T \cdot N_{G}}$ is defined as:
\begin{displaymath}
  \mathbf{a2} \doteq \sum_{g=1}^{N_{G}} \omega_{(t,T_{off})} \boldsymbol{\varrho}_{(t,g)}
\end{displaymath}

Spinning reserve constraints~\eqref{eq:uc:spinning-reserve} are expressed in Eq.~\eqref{eq:uc:m:spinning-reserve} in terms of $\mathcal{A} \mathbf{3} (\cdot)$, which is composed of time matrices $\mathbf{A3}_{t} \in \mathbb{R}^{(2 \cdot T \cdot N_{G}) \times (2 \cdot T \cdot N_{G}) } $,
$\overline{ \mathbf{ P2}} \in \mathbb{R}^{ T \times( T \cdot N_{G} ) }$, and $\mathbf{r} \in \mathbb{R}^{T}$:
\begin{eqnarray*}
  \mathbf{A3}_{t} &\doteq& -\sum_{g=1}^{N_{G}} \boldsymbol{\varrho}_{(t,g)} \boldsymbol{\varrho}_{(t,g)}^{\intercal} \otimes \mathbf{A}\\
  \overline{ \mathbf{ P2}} &\doteq& \sum_{t=1}^{T}  \sum_{g=1}^{N_{G}} \boldsymbol{\xi}_{t} \boldsymbol{\varrho}_{(t,g)}^{\intercal}\\
  \mathbf{r} &\doteq& \left[ SR_{1},\ldots,SR_{T} \right]^{\intercal}
\end{eqnarray*}

Maximum ramp up~\eqref{eq:uc:ramp-up} and down~\eqref{eq:uc:ramp-down} constraints are expressed in Eqs.~\eqref{eq:uc:m:ramp-up} and~\eqref{eq:uc:m:ramp-down}, respectively, in terms of linear mapping $\mathcal{A} \mathbf{4} (\cdot)$, which is composed of time matrices $\mathbf{A4}_{N_{G}(t-1)+g} \in \mathbb{R}^{(2 \cdot T \cdot N_{G}) \times (2 \cdot T \cdot N_{G}) }$:
\begin{displaymath}
  \mathbf{A4}_{N_{G}(t-1)+g} \doteq \left( \boldsymbol{\varrho}_{(t,g)}\boldsymbol{\varrho}_{(t,g)}^{\intercal} - \textbf{1}_{(t>1)} \cdot \boldsymbol{\varrho}_{(t-1,g)}\boldsymbol{\varrho}_{(t-1,g)}^{\intercal} \right) \otimes \mathbf{A}
\end{displaymath}
and $\overline{ \mathbf{ P3}} \in \mathbb{R}^{ (T \cdot N_{B}) \times( T \cdot N_{G} ) }$:
\begin{displaymath}
  \overline{\mathbf{P3}} \doteq \sum_{t=1}^{T} \sum_{g=1}^{N_{G}} \underline{P}_{g} \left( \boldsymbol{\varrho}_{(t,g) }\boldsymbol{\varrho}_{(t,g) } -  \textbf{1}_{(t>1)} \cdot \boldsymbol{\varrho}_{(t-1,g) }\boldsymbol{\varrho}_{(t-1,g) } \right)
\end{displaymath}
as well as $\mathbf{g1}$ and $\mathbf{g2} \in \mathbb{R}^{T}$:
\begin{eqnarray*}
  g1_{\varrho_{4}(t,g)} &\doteq& RU_{g} + \textbf{1}_{(t=1)} \cdot p_{0,g} \\
  g2_{\varrho_{4}(t,g)} &\doteq& RD_{g} - \textbf{1}_{(t=1)} \cdot p_{0,g}
\end{eqnarray*}
for $t=1,\ldots,T$ and $g=1,\ldots,N_{G}$. Where $p_{0,g}$ is active power dispatch at the previous period assumed known.

Power flow constraints~\eqref{eq:uc:flow-limits} are expressed in Eq.~\eqref{eq:uc:m:flow-limits} in terms of linear mapping $\mathcal{Y} \mathbf{3} \left( \cdot \right)$, which is composed of time matrices $\mathbf{Y3}_{r} \in \mathbb{R}^{ (2 \cdot N_{B}) \times( 2 \cdot N_{B})}$,
and $ \overline{\mathbf{f}} \in \mathbb{R}^{2 \cdot T \cdot | \mathcal{L} |}$:
\begin{eqnarray*}
  \mathbf{Y3}_{r} &\doteq& -\mathbf{Y2}_{T \cdot | \mathcal{L} |+r} = \mathbf{Y}_{i,j}\\
  \overline{f}_{r} &\doteq& \overline{f}_{T \cdot | \mathcal{L} |+r} =\overline{F}_{i,j}
\end{eqnarray*}
where $ r=N_{L}(t-1)+l$, for all $ t=1,2,\ldots,T$, and $(i,j) \in \mathcal{L}$.

Voltage constraints~\eqref{eq:uc:voltage} are expressed in Eq.~\eqref{eq:uc:m:voltage} in terms of linear mapping $\mathcal{Y} \mathbf{4} \left( \cdot \right)$, which is composed of time matrices $\mathbf{Y4}_{r} \in \mathbb{R}^{ 2 \cdot N_{B} \times 2 \cdot N_{B}}$,
and $\mathbf{b} \in \mathbb{R}^{2 \cdot T \cdot N_{L}+ 2 \cdot N_{B}}$:
\begin{eqnarray*}
  \mathbf{Y4}_{r} &\doteq& - \mathbf{Y4}_{T \cdot N_{B} +r}\\
   &=& \boldsymbol{\xi}_{i}\boldsymbol{\xi}_{i}^{\intercal} +
  \boldsymbol{\xi}_{N_{B}+i}\boldsymbol{\xi}_{N_{B}+i}^{\intercal}\\
  b_{r} &\doteq& \overline{V}_{i}\\
  b_{ N_{B} \cdot T +r} &\doteq& -\underline{V}_{i}
\end{eqnarray*}
where $r=\left( N_{B}-1\right) t + i $, for all $t=1,2,\ldots,T$, and $i= \{ 1, \ldots , N_{B} \} \backslash i_{slack}$.

%
\subsubsection{Unit Commitment Master Problem}
\label{sec:solution:benders:master}
Problem~\eqref{eq:uc:m:objective}--\eqref{eq:uc:m:binary} is decomposed into a master problem~\eqref{eq:solution:master} and a subproblem~\eqref{eq:solution:subproblem} at every $k$-th iteration, where the latter receives solutions from the former, returning either feasibility or optimality cuts to be added as constraints to the UC master problem, formulated as follows:
\begin{align}
  \min_{\substack{\mathbf{x}, \mathbf{y}, \mathbf{z},w}} \qquad\qquad\quad & w + \mathbf{c}^{\intercal} \mathbf{x} + \mathbf{u}^{\intercal} \mathbf{y} + \mathbf{h}^{\intercal} \mathbf{z} \nonumber\\
  \textrm{s.t.} \qquad\mathbf{Ux} - \mathbf{y} &\leqslant \mathbf{a1} \nonumber\\
  -\mathbf{Ux} - \mathbf{z} &\leqslant - \mathbf{a1} \nonumber\\
  \label{eq:solution:master}%
  \mathbf{M1x} &\leqslant \mathbf{0} \tag{$M^{(k)}$}\\
  \mathbf{M2x} &\leqslant \mathbf{a2} \nonumber\\
  w &\geqslant 0 \nonumber\\
  \mathbf{x}, \mathbf{y}, \mathbf{z} &\in \left\{0,1\right\}^{T\cdot N_{G}} \nonumber
\end{align}

%
\subsubsection{Optimal Power Flow Subproblem}
\label{sec:solution:benders:subproblem}
Subproblem $S^{(k)}$ is parameterized by the solution $\widetilde{\mathbf{x}}^{(k)}$ to master problem $M^{(k)}$, being modified accordingly to return feasibility and optimality cuts. This is accomplished by adding a slack variable $s$ penalized in the objective function by a factor $\sigma\gg 0$. Good numerical stability in the optimization process is obtained in our computational experiments by choosing
\begin{displaymath}
  \sigma \doteq \sum_{g=1}^{N_{G}} \alpha_{g} \cdot \overline{p}_{g}^{2}+\beta_{g} \cdot \overline{p}_{g}+\gamma_{g}
\end{displaymath}
which represents both an upper bound to objective function~\eqref{eq:uc:objective} and a sufficiently big number. If either $s>0$ or its corresponding constraint in the dual problem~\eqref{eq:solution:subproblem-dual} is active in the optimal solution, then~\eqref{eq:solution:subproblem} is infeasible, and the solution ($\mu_{1}^{(k)},\ldots,\mu_{7}^{(k)}$)
to~\eqref{eq:solution:subproblem-dual} is used to construct a feasibility cut~\eqref{eq:feasibility} to be added to~\eqref{eq:solution:master} as a constraint. If otherwise $s=0$ in the optimal solution, then the OPF subproblem~\eqref{eq:solution:subproblem} is feasible and the solution ($\lambda_{1},\ldots,\lambda_{7}$) is used to construct an optimality cut~\eqref{eq:optimality} to be added to~\eqref{eq:solution:master} as a constraint:
\begin{align}
  \min_{\mathbf{P}, \mathbf{V},\mathbf{q}, s} \qquad\qquad\quad \mathbf{C} \bullet \mathbf{P} + \sigma\cdot s \nonumber\\
  \textrm{s.t.} \quad\mathcal{A}\mathbf{1}(\mathbf{P}) + \mathcal{Y} \mathbf{1}(\mathbf{V}) - s\cdot\widehat{\mathbf{e}}_{1}
    &= \mathbf{d}p - \underline{\mathbf{P1}} \widetilde{\mathbf{x}}^{(k)} \nonumber\\
    \quad\mathbf{q} + \mathcal{Y} \mathbf{2}(\mathbf{V}) - s\cdot\widehat{\mathbf{e}}_{2}
    &= \mathbf{d}q - \underline{\mathbf{Q}} \widetilde{\mathbf{x}}^{(k)} \nonumber\\
  \mathcal{A} \mathbf{2}(\mathbf{P}) - s\cdot\widehat{\mathbf{e}}_{3}
    &\leqslant \overline{\mathbf{P1}} \widetilde{\mathbf{x}}^{(k)} \nonumber\\
  \mathbf{q} - s\cdot\widehat{\mathbf{e}}_{4}
    &\leqslant \overline{\mathbf{Q}} \widetilde{\mathbf{x}}^{(k)} \nonumber\\
  \mathcal{A} \mathbf{5}(\mathbf{P})- s\cdot\widehat{\mathbf{e}}_{5}
    &\leqslant \mathbf{r}-\overline{\mathbf{P2}}\widetilde{\mathbf{x}}^{(k)} \nonumber\\
  \label{eq:solution:subproblem}%
  \mathcal{A}\mathbf{4}(\mathbf{P}) - s\cdot\widehat{\mathbf{e}}_{6}
    &\leqslant \mathbf{g1}-\overline{\mathbf{P3}}\widetilde{\mathbf{x}}^{(k)} \tag{$S^{(k)}$}\\
  - \mathcal{A}\mathbf{4}(\mathbf{P}) - s\cdot\widehat{\mathbf{e}}_{7}
    &\leqslant \mathbf{g2}+\overline{\mathbf{P3}}\widetilde{\mathbf{x}}^{(k)} \nonumber\\
  \mathcal{Y}\mathbf{3}(\mathbf{V}) - s\cdot\widehat{\mathbf{e}}_{8}
    &\leqslant \overline{\mathbf{f}} \nonumber\\
  \mathcal{Y}\mathbf{4}(\mathbf{V}) - s\cdot\widehat{\mathbf{e}}_{9}
    &\leqslant \mathbf{b} \nonumber\\
  \mathbf{P}, \mathbf{V}
    &\succeq 0 \nonumber\\
      \mathbf{q}  &\geq 0 \nonumber\\
  s &\geqslant 0 \nonumber
\end{align}

where $\widehat{\mathbf{e}}_{1},\ldots,\widehat{\mathbf{e}}_{9}$ are vectors of ones with appropriate dimension in every equation. In order to obtain a convex subproblem the rank-1 constraint over $\mathbf{V}$ in Eq.~\eqref{eq:uc:m:rank} is relaxed. The corresponding dual of~\eqref{eq:solution:subproblem} is formulated as follows:
\begin{align*}
  \max_{\substack{\boldsymbol{\lambda}_{1},\ldots,\boldsymbol{\lambda}_{9}}} \qquad &
    ( \mathbf{d}p - \underline{\mathbf{P1}} \widetilde{\mathbf{x}}^{(k)} )^{\intercal} \boldsymbol{\lambda}_{1}+( \mathbf{d}q - \underline{\mathbf{Q}} \widetilde{\mathbf{x}}^{(k)} )^{\intercal} \boldsymbol{\lambda}_{2} + ( \overline{\mathbf{P1}} \widetilde{\mathbf{x}}^{(k)} )^{\intercal} \boldsymbol{\lambda}_{3} + ( \overline{\mathbf{Q}} \widetilde{\mathbf{x}}^{(k)} )^{\intercal} \boldsymbol{\lambda}_{4}\\
 & \quad  + ( \mathbf{r} - \overline{\mathbf{P2}}\widetilde{\mathbf{x}}^{(k)} )^{\intercal} \boldsymbol{\lambda}_{5} + ( \mathbf{g1}-\overline{\mathbf{P3}}\widetilde{\mathbf{x}}^{(k)} )^{\intercal} \boldsymbol{\lambda}_{6}
       + ( \mathbf{g2}+\overline{\mathbf{P3}}\widetilde{\mathbf{x}}^{(k)} )^{\intercal} \boldsymbol{\lambda}_{7} +  \overline{\mathbf{f}}^{\intercal} \boldsymbol{\lambda}_{8} + \mathbf{b}^{\intercal} \boldsymbol{\lambda}_{9}
\end{align*}
s.t.
\begin{align}
  \mathcal{A}\mathbf{1}^{\intercal}(\boldsymbol{\lambda}_{1}) + \mathcal{A}\mathbf{2}^{\intercal}(\boldsymbol{\lambda}_{2}) + \mathcal{A}\mathbf{3}^{\intercal}(\boldsymbol{\lambda}_{3})
    + \mathcal{A}\mathbf{4}^{\intercal}( \boldsymbol{\lambda}_{4}) - \mathcal{A}\mathbf{4}^{\intercal}( \boldsymbol{\lambda}_{5})
    &\preceq \mathbf{C4} \nonumber\\
  \label{eq:solution:subproblem-dual}%
  \mathcal{Y}\mathbf{1}^{\intercal}(\boldsymbol{\lambda}_{1}) + \mathcal{Y}\mathbf{2}^{\intercal}(\boldsymbol{\lambda}_{2}) + \mathcal{Y}\mathbf{3}^{\intercal}(\boldsymbol{\lambda}_{6}) + \mathcal{Y}\mathbf{4} ^{\intercal}(\boldsymbol{\lambda}_{7})
    &\preceq \mathbf{0} \tag{$SD^{(k)}$}\\
  - \widehat{\mathbf{e}}_{1}^{\intercal}\boldsymbol{\lambda}_{1} -\cdots - \widehat{\mathbf{e}}_{9}^{\intercal} \boldsymbol{\lambda}_{9}
    &\leqslant \sigma \nonumber\\
  \boldsymbol{\lambda}_{3},\ldots,\boldsymbol{\lambda}_{9}
    &\leqslant \mathbf{0} \nonumber\\
  \boldsymbol{\lambda}_{1} \text{ and } \boldsymbol{\lambda}_{2} &\in \Real^{T\cdot N_{B}} \nonumber
\end{align}

Feasibility cuts are defined below:
\begin{eqnarray}
  \left( \mathbf{d}p - \underline{\mathbf{P1}} \mathbf{x} \right)^{\intercal} \boldsymbol{\mu}_{1}^{(k)}  + \left( \mathbf{d}q - \underline{\mathbf{Q}} \mathbf{x} \right)^{\intercal} \boldsymbol{\mu}_{2}^{(k)}  + \left( \overline{\mathbf{P1}} \mathbf{x} \right)^{\intercal} \boldsymbol{\mu}_{3}^{(k)}+ \left( \overline{\mathbf{Q}} \mathbf{x} \right)^{\intercal} \boldsymbol{\mu}_{4}^{(k)} \notag\\
  + \left( \mathbf{r} - \overline{\mathbf{P2}}\mathbf{x} \right)^{\intercal} \boldsymbol{\mu}_{5}^{(k)} + \left( \mathbf{g1}-\overline{\mathbf{P3}}\mathbf{x} \right)^{\intercal} \boldsymbol{\mu}_{6}^{(k)} +
  \label{eq:feasibility}%
  \left( \mathbf{g2}+\overline{\mathbf{P3}}\mathbf{x} \right)^{\intercal} \boldsymbol{\mu}_{7}^{(k)}+  \overline{\mathbf{f}}^{\intercal} \boldsymbol{\mu}_{8}^{(k)} + \mathbf{b}^{\intercal} \boldsymbol{\mu}_{9}^{(k)} & \geqslant & 0
\end{eqnarray}
as well as optimality cuts:
\begin{eqnarray}
  \left( \mathbf{d}p - \underline{\mathbf{P1}} \mathbf{x} \right)^{\intercal} \boldsymbol{\lambda}_{1}^{(k)}+ \left( \mathbf{d}q - \underline{\mathbf{Q}} \mathbf{x} \right)^{\intercal} \boldsymbol{\mu}_{2}^{(k)}  + \left( \overline{\mathbf{P1}} \mathbf{x} \right)^{\intercal} \boldsymbol{\lambda}_{3}^{(k)} + \left( \overline{\mathbf{Q}} \mathbf{x} \right)^{\intercal} \boldsymbol{\lambda}_{4}^{(k)} \notag\\
  + \left( \mathbf{r} - \overline{\mathbf{P2}}\mathbf{x} \right)^{\intercal} \boldsymbol{\lambda}_{5}^{(k)} + \left( \mathbf{g1}-\overline{\mathbf{P3}}\mathbf{x} \right)^{\intercal} \boldsymbol{\lambda}_{6}^{(k)} +
  \label{eq:optimality}%
  \left( \mathbf{g2}+\overline{\mathbf{P3}}\mathbf{x} \right)^{\intercal} \boldsymbol{\lambda}_{7}^{(k)}+  \overline{\mathbf{f}}^{\intercal} \boldsymbol{\lambda}_{8}^{(k)} + \mathbf{b}^{\intercal} \boldsymbol{\lambda}_{9}^{(k)} &\leqslant& w
\end{eqnarray}

Depending on the nature of the cut generated, it is added to either the set of feasibility cuts:
\begin{displaymath}
  \mathcal{F}_{k} \doteq \left\{ \mathbf{x} \, : \, \textrm{subject to \eqref{eq:feasibility} if } s^{(k-1)} > 0  \right\}
\end{displaymath}
or the set of optimality cuts:
\begin{displaymath}
  \mathcal{O}_{k} \doteq \left\{ \mathbf{x} \, : \, \textrm{subject to \eqref{eq:optimality} if } s^{(k-1)} = 0   \right\}
\end{displaymath}

In OPF problems, zero duality gaps of SDP relaxations are not generally guaranteed, as extensively recorded in the recent literature. In~\cite{6980142} the authors have shown that low order SDP relaxations are sufficient to obtain the global solution of a wide variety of OPF problems. There exist, however, some cases for which SDP relaxations do not guarantee zero duality gaps, as shown in~\cite{6758891} and~\cite{2014arXiv1410.1004K}. Different techniques for dealing with non-zero duality gaps in such problems have been subject of open research. In~\cite{5971792} the authors have shown that minor modifications to line resistances in IEEE systems yield rank-2 solutions that enable further rank-1 solution recovery. It has been observed in~\cite{2014arXiv1410.1004K} that non-convexities in the transmission system can be mitigated by the penalization of apparent losses in order to reduce the difference of angles in the corresponding lines. An analogous approach is presented in~\cite{6980142}, where the authors propose higher order relaxation of the insulting lines. Other alternative is to use SDP relaxations in a branch-and-bound scheme~\cite{6483274}. In this work, a rank reduction procedure is proposed as means to deal with non-zero duality gap of SDP relaxations, as described in Section~\ref{sec:rank}.

%
%
\subsubsection{Modified Master Problem}
\label{sec:solution:benders:modified}
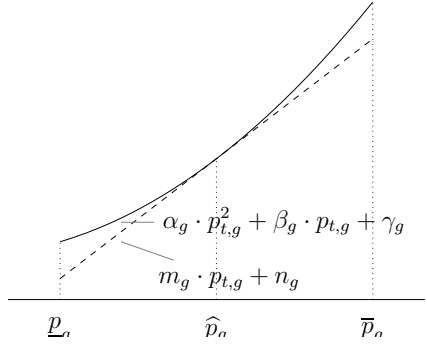
\begin{figure}[!h]
\begin{center}
\begin{tikzpicture}[scale=0.80]
\begin{axis}[axis lines=none,xtick=\empty,ytick=\empty,xmin=0,xmax=80,ymin=-1900,ymax=6500]
\addplot[ black, domain=10:70]{x^2+17.7*x+137};
\draw[dotted] (axis cs:10,-1000) -- (axis cs:10,414);
\draw[dotted] (axis cs:40,-1000) -- (axis cs:40,2445);
\draw[dotted] (axis cs:70,-1000) -- (axis cs:70,6276);
\draw[black] (axis cs:0,-1000) -- (axis cs:80,-1000);
\addplot[ black, dashed, domain=10:70]{97.7*x-1600+137};
\addplot[] coordinates {(20,891)} node[pin=0:{$\alpha_{g} \cdot  p_{t,g}^{2}+\beta_{g} \cdot p_{t,g}+\gamma_{g} $}]{} ;
\addplot[] coordinates {(20,491 )} node[pin=-30:{$m_{g}\cdot p_{t,g}+n_{g}$}]{} ;
\node[label={-90:{$\underline{p}_{g}$}}] at (axis cs:10,-1000) {};
\node[label={-90:{$\widehat{p}_{g}$}}] at (axis cs:40,-1000) {};
\node[label={-90:{$\overline{p}_{g}$}}] at (axis cs:70,-1000) {};
\end{axis}
\end{tikzpicture}
\end{center}
 \caption{Linear lower bound}\label{fig:llowerb}
\end{figure}

A linearization of~\eqref{eq:solution:subproblem} is performed into the modified master problem~\eqref{eq:solution:modified-master} as a means to represent the subproblem into the master problem, thus constraining integer solutions, and consequently improving convergence of the overall problem, especially in the initial iterations, where, if we don't use this linearization in the master problem, it will result in many infeasible unit commitment solutions to the subproblems. At some point of the iterations it's expected that this linearization will be less tight than the Benders' cuts and could also be neglected. Therefore, variables $p_{t,g}$ are introduced in the master problem, as well as a corresponding relaxed model representing active power generation constraints from the OPF subproblem, consisting of relaxed power balance equations:
\begin{equation}
  \label{eq:active-demand-supply}%
\sum_{g=1}^{N_{G}}p_{t,g}-l_{t}^{(k-1)} \geqslant \sum_{i=1}^{N_{B}}dp_{t,i}
\end{equation}
\begin{equation}
  \label{eq:reactive-demand-supply}%
\sum_{g=1}^{N_{G}}q_{t,g} \geqslant \sum_{i=1}^{N_{B}}dq_{t,i}
\end{equation}
for $t=1,2,\ldots,T$, where $l_{t}^{(k-1)}$, representing power losses updated from ($S^{k-1}$), is assumed as a ratio of total load at the beginning of each iteration. At optimality, Eq.~\eqref{eq:active-demand-supply} and Eq.~\eqref{eq:reactive-demand-supply}  are active. Also, spinning reserve constraints
\begin{equation}
  \label{eq:spinningMax}%
   \sum_{i=1}^{N_{B}}dp_{t,i} + l_{t}^{(k-1)} + SR_{t} \leqslant \sum_{g=1}^{N_{G}} \overline{p}_{g}\cdot x_{t,g}
\end{equation}
the equivalent for reactive generation:
\begin{equation}
  \label{eq:reactiveMax}%
   \sum_{i=1}^{N_{B}}dq_{t,i} + l_{t}^{(k-1)} \leqslant \sum_{g=1}^{N_{G}} \overline{q}_{g}\cdot x_{t,g}
\end{equation}
for $t=1,2,\ldots,T$, minimum and maximum active and reactive power generation
\begin{equation}
  \label{eq:limits}%
  \underline{p}_{g} \cdot x_{t,g} \leqslant p_{t,g} \leqslant \overline{p}_{g} \cdot x_{t,g}
\end{equation}
\begin{equation}
  \label{eq:limitsq}%
  \underline{q}_{g} \cdot x_{t,g} \leqslant q_{t,g} \leqslant \overline{q}_{g} \cdot x_{t,g}
\end{equation}
ramp up
\begin{equation}
  \label{eq:rampup}%
  p_{t,g}-p_{t-1,g} \leqslant RU_{g}
\end{equation}
and, finally, ramp down constraints
\begin{equation}
  \label{eq:rampdown}%
  p_{t-1,g}-p_{t,g} \leqslant RD_{g}
\end{equation}
for $t=1,2,\ldots,T$ and $g=1,2,\ldots,N_{G}$. Additionally, objective function~\eqref{eq:uc:objective} is approximated by a linear function as illustrated in Fig.~\ref{fig:llowerb}. This linear approximation represents a lower bound to the quadratic fuel costs, and is calculated by using the mean value theorem of convex functions:
\begin{equation*}
  \left( 2 \alpha_{g} \cdot \widehat{p}_{g} +\beta \right) \left( p_{t,g}- \widehat{p}_{g} \right) +
  \alpha_{g} \cdot \widehat{p}_{g}^{2}+\beta_{g} \cdot \widehat{p}_{g} + \gamma_{g}  \leqslant  \alpha_{g} \cdot  p_{t,g}^{2}+\beta_{g} \cdot p_{t,g} +\gamma_{g}
\end{equation*}

Furthermore, without any loss of generality, we choose:
\begin{displaymath}
  \widehat{p}_{g} \doteq \frac{\overline{p}_{g}-\underline{p}_{g}}{2}
\end{displaymath}
and finally obtain a linear lower bound of the cost function represented by $w+\mathbf{c}^{\intercal} \mathbf{x}$:
\begin{equation}
  \label{eq:llowerb}%
  \mathbf{m}^{\intercal} \mathbf{p} + \mathbf{n}^{\intercal} \mathbf{x} \leqslant  w + \mathbf{c}^{\intercal} \mathbf{x}
\end{equation}
where:
\begin{align*}
  \mathbf{m} &\doteq \left[m_{1},\ldots,m_{N_{G}},\ldots,m_{1},\ldots,m_{N_{G}} \right]^{\intercal} \, \in \mathbb{R}^{T\cdot N_{G}}\\
  \mathbf{n} &\doteq \left[n_{1},\ldots,n_{N_{G}},\ldots,n_{1},\ldots,n_{N_{G}}  \right]^{\intercal} \, \in \mathbb{R}^{T\cdot N_{G}}
\end{align*}
and
\begin{align*}
  m_{g} &= 2 \alpha_{g}\cdot \widehat{p}_{g} +\beta_{g}\\
  n_{g} &=- \alpha_{g}\cdot \widehat{p}_{g}^2 +\gamma_{g}
\end{align*}
for $t=1,2,\ldots,T$, and $g=1,2,\ldots,N_{G}$.

Finally, the modified master problem is formulated below:\begin{equation}
  \label{eq:solution:modified-master}%
  \begin{array}{crclcl}
    \displaystyle \min_{ \substack{ \mathbf{x},\mathbf{y},\mathbf{z}, \\ \mathbf{p},\mathbf{q},w}}  & \multicolumn{5}{l}{ w + \mathbf{c}^{\intercal} \mathbf{x} +\mathbf{u}^{\intercal} \mathbf{y}+ \mathbf{h}^{\intercal} \mathbf{z} } \\
    \textrm{s. t. } & \mathbf{Ux} -\mathbf{y} &\leqslant&\mathbf{a1} & & \\
     & -\mathbf{Ux} - \mathbf{z}&\leqslant&-\mathbf{a1} & & \\
     & \mathbf{M1x} &\leqslant&\mathbf{0} & & \\ \tag{$MM^{(k)}$}
     & \mathbf{M2x} &\leqslant&\mathbf{a2} & & \\
     & \mathbf{H} \mathbf{p} + \mathbf{G}\mathbf{x} & \leqslant & \bm{\pi}\\
     & \widetilde{\mathbf{H}} \mathbf{q} + \widetilde{\mathbf{G}}\mathbf{x} & \leqslant & \widetilde{\bm{\pi}}\\
     & \mathbf{m}^{\intercal} \mathbf{p} + \mathbf{n}^{\intercal} \mathbf{x} & \leqslant & w + \mathbf{c}^{\intercal}\mathbf{x}\\
     & \mathbf{p}, \; \mathbf{q} & \geqslant & 0\\
     & \mathbf{x} & \in & \cap_{i=1}^{k} \{ \mathcal{F}_{i} \cup \mathcal{O}_{i} \}\\
     & \mathbf{x},\mathbf{y},\mathbf{z}  & \in & \left\{ 0,1 \right \}^{T \cdot N_{G}}\\
  \end{array}
\end{equation}

where $\mathbf{H} \mathbf{p} + \mathbf{G}\mathbf{x} \leq \bm{\pi}$ represents constraints \eqref{eq:active-demand-supply}, \eqref{eq:spinningMax}, \eqref{eq:limits}, \eqref{eq:rampup}, and \eqref{eq:rampdown}, $\widetilde{\mathbf{H}} \mathbf{q} + \widetilde{\mathbf{G}}\mathbf{x}  \leqslant  \widetilde{\bm{\pi}}$ represents constraints \eqref{eq:reactive-demand-supply}, \eqref{eq:reactiveMax}, \eqref{eq:limitsq}, and feasibility~$\mathcal{F}_{i}$ and optimality~$\mathcal{O}_{i}$ cuts are added iteratively. Note that $w\geqslant 0$ is incurred from the cost function lower bound in Eq.~\eqref{eq:llowerb}.

%
%
\subsection{Rank Reduction Procedure}
\label{sec:rank}
As demonstrated in Section~\ref{sec:problem:constraints:rank}, the rank of $\mathbf{P}^{*}$ is guaranteed to be one at every iteration. A rank reduction procedure~\cite{1997geometry} is necessary, however, such that a perturbation matrix $\mathbf{D}_{r}$ in the rank reduction iteration $r$ is added to the voltage matrix $\mathbf{V}^{*}$ that solves problem~\eqref{eq:uc:m:objective}--\eqref{eq:uc:m:binary} at the optimum of the Bender's algorithm. The rank reduction procedure listed in Algorithm~\ref{alg:alg1} maintains feasibility as explained in the remainder of this section. Moreover, because $\mathbf{V}$ does not appear in the objective function, new solutions resulting from the rank reduction procedure are still optimal. A central idea to this procedure is to calculate $\mathbf{D}_{r}$ such that at every iteration, $\mathbf{V}^{*}_{r}$ does not violate any constraints, has its rank reduced, and be positive semidefinite.
\begin{algorithm}
\begin{center}
\begin{algorithmic}[1]
\REQUIRE $\mathbf{V}^{*}$,$ \{ \mathbf{Y}_{i} \forall i \in \{1,\ldots,N_{B} \}\}$, $\{ \mathbf{Y}_{i,j} \forall (i,j) \in \mathcal{L} \}$ and  $\{\mathbf{E}_{l} \forall l \in \{1,\ldots,2N_{B}  \} \}$
\STATE $\mathbf{V}^{*}_{1}=\mathbf{V}^{*}$
\FOR{$r\leftarrow 1,\ldots,\left( \text{rank}(\mathbf{V}^{*})+\varsigma-1 \right) $}
\STATE Obtain Cholesky decomposition $\mathbf{V}^{*}_{r}=\mathbf{R}_{r}\mathbf{R}_{r}^{\intercal}$
\STATE Initialize $\mathbf{S}_{r} \leftarrow \mathbf{0}_{2N_{B} \times 2N_{B} }$
\FOR{$i\leftarrow 1,\ldots,N_{B}$ }
\STATE  $\mathbf{S}_{r} \leftarrow \mathbf{S}_{r} + \left| \mathbf{R}_{k}^{\intercal}\mathbf{Y}_{i} \mathbf{R}_{r} \right|$
\ENDFOR
\FOR{$(i,j) \leftarrow \mathcal{L} $ }
\STATE  $\mathbf{S}_{r} \leftarrow \mathbf{S}_{r}+ \left| \mathbf{R}_{r}^{\intercal} \mathbf{Y}_{i,j} \mathbf{R}_{r}\right|$
\ENDFOR
\FOR{$l \leftarrow 1,\ldots,2N_{B}$ }
\STATE  $\mathbf{S}_{r} \leftarrow \mathbf{S}_{r}+ \left| \mathbf{R}_{r}^{\intercal} \mathbf{E}_{l} \mathbf{R}_{r} \right|$
\ENDFOR
\FORALL{$S_{i,j} \in \mathbf{S}_{r} \, | \, \mathbf{Z}_{r} \in \mathbb{R}^{2N_{B}\times2N_{B}}$}
\STATE $Z_{i,j} \leftarrow \textbf{1}_{(S_{i,j}= 0)}$
\ENDFOR
\IF{$\mathbf{Z}_{r} = \mathbf{0}$}
\STATE \textbf{break}
\ENDIF
\STATE $\phi_{z} \leftarrow 2 \cdot \textbf{1}_{(\mathbf{z}\succeq \mathbf{0})}-1$
\STATE $\mathbf{D}_{r} \leftarrow -\phi_{z} \cdot \mathbf{R}_{r}\mathbf{Z}_{r}\mathbf{R}_{r}^{\intercal}$
\STATE $\omega_{r} \leftarrow \min \left\{ \, \left(\phi_{z} \cdot \lambda_{z}\right)^{-1} \, : \, z = 1, 2, \ldots, N_{B} \, \right\}$
\STATE $\mathbf{V}^{*}_{r+1} \leftarrow \mathbf{V}^{*}_{r} + \omega_{r}  \cdot \mathbf{D}_{r}$
\IF{$\text{rank}(\mathbf{V}^{*}_{r+1}) = 1$}
\STATE \textbf{break}
\ENDIF
\ENDFOR
\end{algorithmic}
\end{center}
\caption{Rank reduction procedure.}\label{alg:alg1}
\end{algorithm}

An upper bound on the rank of the final optimal solution to subproblem~\eqref{eq:solution:subproblem}  can be expressed as a function of the number of problem constraints, as shown in~\cite{barvinok2002course}:

\begin{equation}
  \label{eq:rankmax}%
  \text{rank}(\mathbf{V}^{*}_{r}) \leqslant \text{rank}(\mathbf{V}^{*}_{1})
   \leqslant \left\lfloor \frac{\sqrt{8\left( N_{B}+2|\mathcal{L} |+2N_{B} \right)+1}-1}{2} \right \rfloor
\end{equation}
such that $\mathbf{V}^{*}_{r}$ is updated as follows:
\begin{equation}
  \label{eq:rankup}%
  \mathbf{V}^{*}_{r+1} \leftarrow \mathbf{V}^{*}_{r} +\omega_{r}\cdot \mathbf{D}_{r}
\end{equation}
where $\omega_{r}$ is a scalar parameter that guarantees that every update to $\mathbf{V}^{*}_{r}$ remains positive semidefinite. It follows that the perturbation matrix $\mathbf{D}_{r}$ must be subject to the following constraints:
\begin{align}
  \label{eq:feas-Yi}%
  \mathbf{D}_{r} \bullet \mathbf{Y}_{i} & = 0, \, i=1,\ldots,N_{B}  \\
  \label{eq:feas-Yij}%
  \mathbf{D}_{r} \bullet \mathbf{Y}_{i,j} & = 0, \, \left(i,j\right) \in \mathcal{L} \\
  \label{eq:feas-El}%
  \mathbf{D}_{r} \bullet \mathbf{E}_{l} & = 0, \, l = 1,2,\ldots,2N_{B}
\end{align}
where $\mathbf{D}_{r}$ has the following form:
\begin{equation}
  \label{eq:decomp}%
  \mathbf{D}_{r} = \mathbf{R}_{r} \mathbf{Z}_{r} \mathbf{R}_{r}^{\intercal}
\end{equation}
where $\mathbf{R}_{r}$ is obtained from the Cholesky decomposition of $\mathbf{V}^{*}_{r} = \mathbf{R}_{r}\mathbf{R}_{r}^{\intercal}$. By the distributive property of the Frobenious product we obtain the following relationship:
\begin{align}
  \notag \mathbf{D}_{r}\bullet\mathbf{Y}_{i} & = \mathbf{R}_{r} \mathbf{Z}_{r} \mathbf{R}_{r}^{\intercal} \bullet \mathbf{Y}_{i}
  \notag  = \mathbf{Z}_{r} \bullet \mathbf{R}_{r}^{\intercal}  \mathbf{Y}_{i}  \mathbf{R}_{r}
  \notag  = 0\\
  \label{eq:distributive-property}%
   & \rightarrow \mathbf{Z}_{r} \left( \mathbf{R}_{r}^{\intercal}  \mathbf{Y}_{i}  \mathbf{R}_{r} \right) = \mathbf{0}
\end{align}

In order to satisfy~\eqref{eq:distributive-property}, $\mathbf{Z}_{r}$ is constructed as a $0$-$1$ matrix such that there is a one-valued element for every zero-valued element in the corresponding coordinate of matrices $\left( \mathbf{R}_{r}^{\intercal}  \mathbf{Y}_{i}  \mathbf{R}_{r} \right)$, $\left( \mathbf{R}_{r}^{\intercal}  \mathbf{Y}_{i,j}  \mathbf{R}_{r} \right)$, and $\left( \mathbf{R}_{r}^{\intercal}  \mathbf{E}_{l}  \mathbf{R}_{r} \right)$. This guarantees that~\eqref{eq:feas-Yi}, \eqref{eq:feas-Yij}, and~\eqref{eq:feas-El} are satisfied.

Substituting~\eqref{eq:decomp} in~\eqref{eq:rankup}, we obtain:
\begin{equation}
  \label{eq:rankred}%
  \mathbf{V}^{*}_{r} + \omega_{r} \cdot \mathbf{D}_{r} = \mathbf{R}_{r} \left( \mathbf{I} - \omega_{r}\mathbf{Z}_{r} \right) \mathbf{R}_{r}^{\intercal}
\end{equation}
From~\eqref{eq:rankred} we conclude that the necessary condition for positive semidefiniteness of the new solution is given by:
\begin{displaymath}
  \left( \mathbf{V}^{*}_{r} + \omega_{r}\cdot \mathbf{D}_{r} \right)   \succeq \mathbf{0} \Leftrightarrow \left( \mathbf{I} - \omega_{r}\mathbf{Z}_{r} \right) \succeq \mathbf{0}
\end{displaymath}

If $\mathbf{Z}_{r} \preceq \mathbf{0}$, then a multiplication by $-1$ follows. This is characterized by the value assumed by parameter $\phi_{z}$. If $\mathbf{Z}_{r}$ is indefinite, then the algorithm is not able to reduce its rank at the $r$-th iteration. However it is still possible to reduce its rank in the following iterations since $\mathbf{Z}_{r}\neq\mathbf{Z}_{r+1}$ because of the added perturbations, and also because the SDP formulation of the network constraints guarantees that at least one rank-1 solution exists. In that sense, additional $\varsigma$ iterations are provisioned. Parameter $\omega_{r}$, on the other hand, is a scale of perturbation matrix $\mathbf{D}_{r}$ that induces rank reduction of $\mathbf{V}^{*}_{r}$ in at least one degree by reducing the rank of $\left( \mathbf{I}- \omega_{r}\mathbf{Z}_{r} \right)$. This is achieved by:
\begin{displaymath}
  \omega_{r} = \min \left\{ \, \left(\phi_{z} \cdot\lambda_{z}\right)^{-1} \, : \, z = 1, 2, \ldots, \text{rank}(\mathbf{V}_{r}^{*}) \, \right\}
\end{displaymath}
where $\lambda_{z}$ is an eigenvalue of $\mathbf{Z}_{r}$.

\begin{remark}
	If $\text{rank}(\mathbf{V}_{r}^{*})>1$ then feasibility of the solution must be recovered. This can be achieved, without any guarantees of global optimality, by means of a local resolution, either approximately by, e.g. a convex-concave procedure~\cite{Lipp2016}, or locally, using a method to solve the AC OPF with strong local convergence properties, e.g. an infeasible primal-dual interior-point method with line search filter~\cite{Wachter2006}.
\end{remark}

\section{Numerical Results}
\label{sec:results}
The Benders'-like decomposition algorithm described in Section~\ref{sec:solution:benders} was implemented in MATLAB$^\circledR$ 8.1 with solution of the mixed-integer master problem by the IBM$^\circledR$ CPLEX$^\circledR$ 12.6 solver using the branch-and-cut algorithm, and solution of the OPF subproblem by SDP relaxation using SDPA 7.3.8.

\begin{table}[!t]
\centering
\caption{Generator Data for the 3-GEN case.}\label{tab:3gendata}
\pgfplotstableset{every head row/.style={before row=\toprule, after row=\midrule}, every last row/.style={after row=\bottomrule},}
\begin{filecontents*}{3-GENTgendata.dat}
Parameter	G1	G2	G3	
alpha	0.0004	0.0010	0.0050	
beta	13.70	40.00	17.70	
gamma	177.00	130.00	137.00	
Pmax	210.00	100.00	70.00
Pmin	100.00	10.00	10.00
Qmax	210.00	100.00	70.00	
Qmin	-210.0	-100.0	-70.0		
SU	100.00	200.00	50.00	
SD	50.00	100.00	50.00	
Pp	150.00	0.00	15.00	
UR	55.00	50.00	15.00	
DR	55.00	50.00	15.00	
To	2	-1	1	
Ton	4	3	2	
Toff	4	2	2	
\end{filecontents*}
\pgfplotstabletypeset[ 1000 sep={,}, string replace={Pmax}{$\overline{P}$ $\left[\si{\mega \watt } \right]$},string replace={Pmin}{$\underline{P}$ $\left[\si{\mega \watt } \right]$},string replace={Qmax}{$\overline{Q}$ $\left[\si{\mega VAR } \right]$},string replace={Qmin}{$\underline{Q}$ $\left[\si{\mega VAR } \right]$},string replace={To}{$T_{0}$ $\left[\si{\hour} \right]$},string replace={Ton}{$T_{on}$ $\left[\si{\hour} \right]$},string replace={Toff}{$T_{off}$ $\left[\si{\hour } \right]$},string replace={Pp}{$p_{0}$ $\left[\si{\mega \watt} \right]$},string replace={SU}{$SU$ $\left[\$ \right]$},string replace={SD}{$SD$ $\left[\$ \right]$},string replace={UR}{$UR$ $\left[\si{\mega \watt} \right]$},string replace={DR}{$DR$ $\left[\si{\mega \watt} \right]$},string replace={alpha}{$\alpha$ $\left[\si{\$\per \mega \watt^2} \right]$},string replace={beta}{$\beta$ $\left[\si{\$\per \mega \watt} \right]$}, string replace={gamma}{$\gamma$ $\left[\$\right]$},columns/Parameter/.style={string type,column type=l},columns/info/.style={fixed,fixed zerofill,precision=1,showpos, column type=r,},columns={ [index]0,
[index]1,[index]2,[index]3}]{3-GENTgendata.dat}
\end{table}

In order to validate the methodology we have constructed two numerical cases: 3-GEN and IEEE-118. The purpose of the small-scale 3-GEN case is to didactically serve to showcase the proposed solution to the UC problem. It has a 6-bus network with limits on power flow, and 3 generators with different cost functions, as listed on Table~\ref{tab:3gendata}. In this case generator G1 is less costly and has higher capacity, whereas generators G2 and G3 present equivalent costs depending on the range of operation. On the other hand, the purpose of the case study based on the IEEE 118-bus dataset\footnote{Problem data can be obtained from the following file available online: \url{http://motor.ece.iit.edu/Data/SCUC_118.xls} and also by request to the authors}. is to test the effectiveness of the proposed methodology in a larger-scale problem. In order to speed up the convergence of the algorithm, we adopted a strategy that assumes initial estimated power losses to be between 5\% and 10\% of total system load.

\begin{figure}[!h]
  \begin{center}
    \begin{tiny}
      \begin{tikzpicture}
        \begin{axis}[xlabel=Iterations,ylabel=Operational Cost (\$),width=8.5cm, height=4cm,
            axis x line*=bottom,
    axis y line*=left,
    xmin=1,xmax=6,ymin=88000,ymax=115000,legend style={at={(0.8,0.9)}, anchor=north,legend columns=-1}]
           \pgfplotstableread{UBLB6bus.txt}\datatable;
          \addplot [dashed, color=gray,const plot , mark = ] table[x index=0 , y index=1] from \datatable;
          \addlegendentry{LB};
          \addplot [color=black,const plot , mark = ] table[x index=0 , y index=2] from \datatable;
          \addlegendentry{UB};
        \end{axis}
      \end{tikzpicture}
    \end{tiny}
  \end{center}
  \caption{Upper and lower bound convergence for the 3-GEN case.}
  \label{fig:UBLB3GENR1}
\end{figure}
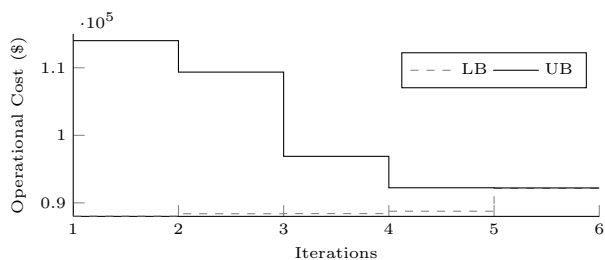

In Fig.~\ref{fig:UBLB3GENR1} we can observe the evolution of both the upper and lower bounds during the algorithm iterations while solving the 3-GEN case. As expected, the gap between the two bounds is very large in the first iteration. However, in only 6 iterations upper and lower bounds converge to an infinitesimal gap with 3 feasibility cuts and 3 optimality cuts.

\begin{table}[!h]
    \caption{Iterations until convergence with and without modified master problem.}
  \centering
  \begin{threeparttable}
    \begin{tabular}{lcc}
      \toprule
      Case study & $MM^{k}$ & $M^{k}$ \\
      \midrule
      3-GEN & 6 & 13 \\
      IEEE-118 & 6 & $\infty$ \\
      \bottomrule
    \end{tabular}
  \end{threeparttable}
  \label{tab:MMP}%
\end{table}

\pgfplotstableread{./ucsix.dat}\datatable
\pgfplotsset{width=12cm}

\begin{figure}[!t]
  \centering
  \begin{tiny}
    \begin{tikzpicture}
      \begin{axis}[domain=1:24,
          view={40}{20},
          axis x line*=bottom,
          axis y line*=left,
          axis z line=none,
          unit vector ratio*=1 3 1,
          xtick={1,...,244}, ytick={1,...,6}, ztick=\empty,
          x tick style = transparent,
          y tick style = transparent,
          x tick label style={rotate=45,anchor=east},
          yticklabels={1,...,6},
          xlabel=Hours,
          x label style={rotate=-15,anchor=east},
          ylabel=Iterations,
          y label style={rotate=20,anchor=west},
          clip=false,
        ]

        \def\minuscurve{0}
        \pgfplotsinvokeforeach{1,...,24}{
          \draw [gray!20] (axis cs:#1,1,0) -- (axis cs:#1,6,0);
        }

        \addplot3[fill=lightgray,mark=none,gray!20,domain=1:24] coordinates {
          (9,1,0)
          (10,1,0)
          (10,6,0)
          (9,6,0)
          (9,1,0)};

        \addplot3[fill=lightgray,mark=none,gray!20,domain=1:24] coordinates {
          (22,1,0)
          (23,1,0)
          (23,6,0)
          (22,6,0)
          (22,1,0)};

        \def\sumcurve{0}
        \pgfplotsinvokeforeach{1,...,6}{
          \draw [gray!20] (axis cs:1,#1,0) -- (axis cs:24,#1,0);
          \addplot3[ycomb,mark=*,mark size=1.0pt,domain=1:24] table [x=T, y expr=#1, z=UC#1] \datatable;
        }
      \end{axis}
    \end{tikzpicture}
  \end{tiny}
  \caption{Iterative progress of the unit commitment schedule of generator G2 for the 3-GEN case.}
    \label{fig:unitx}%
\end{figure}
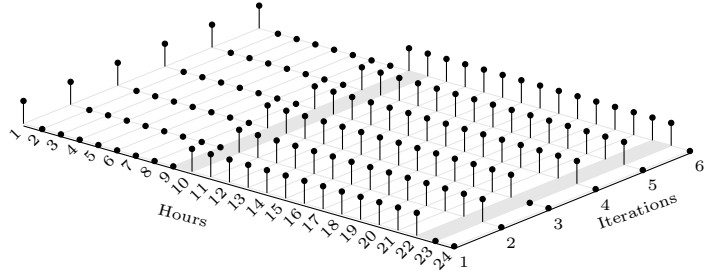

A characteristic of Benders'-like decomposition methods is that the solutions of the integer master problem are evaluated as they serve as parameters of the subproblem. This interaction between the master and subproblem explains the behaviour observed in Fig.~\ref{fig:unitx}, which illustrates how the unit commitment schedule of generator G2 is modified over the iterations as it adapts to the inclusion of power system constraints---e.g. in this case, feasibility cuts associated with spinning reserve and ramp rate constraints.

\begin{figure}
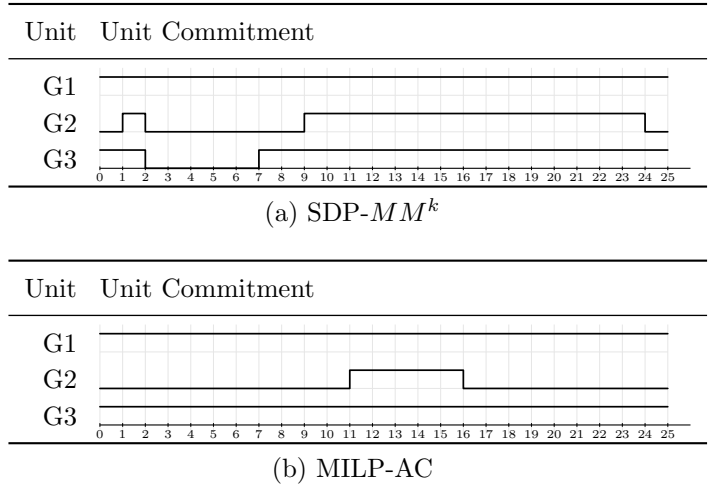
%
	\centering
	\subfloat[][SDP-$MM^k$]{
      \begin{footnotesize}
        \begin{tikztimingtable}[timing/slope=0, font=\rmfamily,timing/xunit=8.5 ]
          G1 &  25{H} \\
          G2 & 1L1H7L15{H}L  \\
          G3 &  2H5L18{H} \\
          \begin{extracode}
            \begin{background}
              \vertlines[help lines,gray!20]{}
              \horlines[help lines,gray!20]{}
              \draw [-,>=latex] (0,-\nrows-1) -- (\twidth+1,-\nrows-1);
              \foreach \n in {0,1,...,\twidth}
              \draw (\n,-\nrows-1+.1) -- +(0,-.2)
              node [below,inner sep=1pt] {\scalebox{.75}{\tiny\n}};
            \end{background}
            \tableheader[]{ Unit }{Unit Commitment}
            \tablerules
          \end{extracode}
        \end{tikztimingtable}
      \end{footnotesize}%
	\label{fig:unitc-sdp}}
	
	\subfloat[][MILP-AC]{
      \begin{footnotesize}
        \begin{tikztimingtable}[timing/slope=0, font=\rmfamily,timing/xunit=8.5 ]
          G1 & 25{H}\\
          G2 & 11L5H9L\\
          G3 & 25{H}\\
          \begin{extracode}
            \begin{background}
              \vertlines[help lines,gray!20]{}
              \horlines[help lines,gray!20]{}
              \draw [-,>=latex] (0,-\nrows-1) -- (\twidth+1,-\nrows-1);
              \foreach \n in {0,1,...,\twidth}
              \draw (\n,-\nrows-1+.1) -- +(0,-.2)
              node [below,inner sep=1pt] {\scalebox{.75}{\tiny\n}};
            \end{background}
            \tableheader[]{ Unit }{Unit Commitment}
            \tablerules
          \end{extracode}
        \end{tikztimingtable}
      \end{footnotesize}%
	\label{fig:unitc-milp}}
	\caption{Optimal unit commitment schedule for the 3-GEN case.}
	\label{fig:unitc}
\end{figure}

\begin{figure}
  \centering
  \subfloat[][SDP-$MM^k$]{
  \begin{footnotesize}
    \begin{tikzpicture}
      \begin{axis}[xlabel=Hours,ylabel=Power (MW),xtick=data,ytick={10,50,...,290},width=12cm, height=4cm ,xmin=1,xmax=24,ymin=0,ymax=290,legend style={at={(0.5,-0.45)}, anchor=north,legend columns=-1}]
        \pgfplotstableread{6busresult.dat}\datatable;
        \addplot [dashed, mark = ] table[y=Pg1] from \datatable ;
        \addlegendentry{G1};
        \addplot [color = gray, mark = ] table[y=Pg2] from \datatable ;
        \addlegendentry{G2};
        \addplot [dashed,color = gray, mark = ] table[y=Pg3] from \datatable ;
        \addlegendentry{G3};
        \addplot [color = black, mark = ] table[y expr=\thisrowno{2} + \thisrowno{5} + \thisrowno{8}] from \datatable;
        \addlegendentry{Demand};
      \end{axis}
    \end{tikzpicture}
  \end{footnotesize}%
  \label{fig:gend-sdp}}
  
  \subfloat[][MILP-AC]{
  \begin{footnotesize}
    \begin{tikzpicture}
      \begin{axis}[xlabel=Hours,ylabel=Power (MW),xtick=data,ytick={10,50,...,290},width=12cm, height=4cm ,xmin=1,xmax=24,ymin=0,ymax=290,legend style={at={(0.5,-0.5)}, anchor=north,legend columns=-1}]
        \pgfplotstableread{6busresult-milp.dat}\datatable;
        \addplot [dashed, mark = ] table[y=Pg1] from \datatable ;
        \addlegendentry{G1};
        \addplot [color = gray, mark = ] table[y=Pg2] from \datatable ;
        \addlegendentry{G2};
        \addplot [dashed,color = gray, mark = ] table[y=Pg3] from \datatable ;
        \addlegendentry{G3};
        \addplot [color = black, mark = ] table[y expr=\thisrowno{2} + \thisrowno{5} + \thisrowno{8}] from \datatable;
        \addlegendentry{Demand};
      \end{axis}
    \end{tikzpicture}
  \end{footnotesize}%
  \label{fig:gend-milp}}
  \caption{Optimal hourly active power generation for the 3-GEN case.}
  \label{fig:gend}
\end{figure}
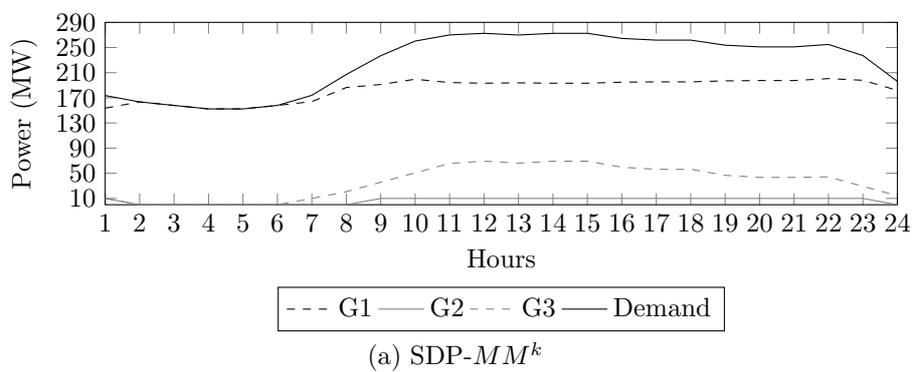
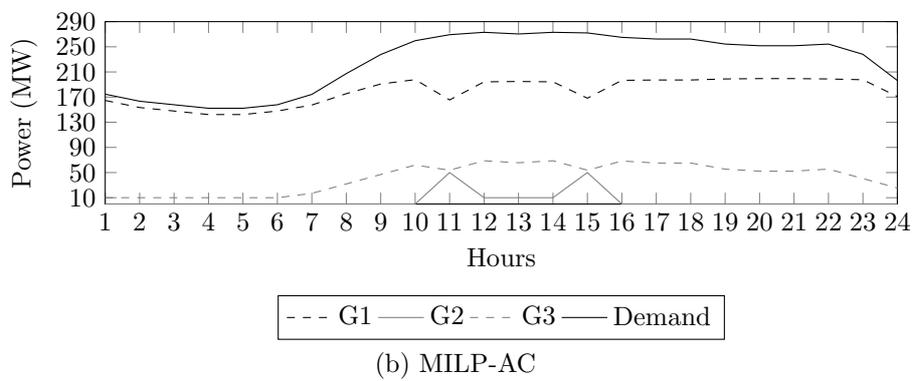

Unit commitment and power generation results are illustrated in Fig.~\ref{fig:unitc-sdp} and Fig.~\ref{fig:gend-sdp}, respectively. It can be observed in Table~\ref{tab:3gendata} that generator G1 is the least costly, resulting in it being dispatched at all times as illustrated in Fig.~\ref{fig:unitc}, whereas more costly generators G2 and G3 complement generation to meet the load demand and spinning reserve requirements. This observed behaviour in the final solution is a numerical evidence of consistency towards solution optimality and feasibility. In order to illustrate the effectiveness of the proposed methodology, it is compared with a simple heuristic deemed herein as MILP-AC. In this heuristic, the unit commitment problem is solved for linearized network constraints with DC power flow, and piecewise-linear costs as a mixed-integer linear program whose unit commitment results are shown in Fig.~\ref{fig:unitc-milp}. Upon optimality, a power factor of .85 is used to calculate active and reactive power generation set points from the optimal DC power flow. These set points are used to initialize the AC power flow for the final generation dispatches as shown in Fig.~\ref{fig:gend-milp}. The MILP-AC heuristic was implemented in the Python programming language using PyPSA~\cite{PyPSA} and Pyomo~\cite{hart2011pyomo} with IBM$^\circledR$ CPLEX$^\circledR$. Pyomo was used to define the piecewise-linear cost functions, and spinning reserve and initial ramping constraints not readily available in PyPSA version 0.17 at the time of writing. In the 3-GEN case, both the Benders' $MM^k$ decomposition with SDP relaxation of AC power flow constraints and the MILP-AC heuristic yielded UC and active power feasible solutions, with a cost difference of 0,52\% lower for MILP-AC at \$93,404.52.

In the IEEE-118 case, generation unit configurations similar to those of the 3-GEN case were used. Convergence of the larger-scale case is also achieved in only six iterations as illustrated in Table~\ref{tab:ieee118sol1}. In order to test the proposed methodology for its effectiveness in solving the AC optimal power flow, we reduced the limits of some power lines. This additional restraining is dealt with effectively by the convex SDP relaxation. As a comparison, in this modified IEEE 118-bus case the MILP-AC solution is infeasible for active power generation after the ensuing AC power flow with a 14\% cost increase at \$1,958,123.15.

The low number of iterations required for convergence of our Benders'-like decomposition algorithm is a result of the use of the modified master problem~\eqref{eq:solution:modified-master}, as evidenced by additional numerical experiments. We verified that convergence was reached in 13 iterations for the same 3-GEN case study if modification of the master problem is not present as described in Section~\ref{sec:solution:benders:modified}, as illustrated in Table~\ref{tab:MMP}. It can be observed that convergence of the IEEE-118 case is not even possible without the use of~\eqref{eq:solution:modified-master}. It is important to note that in these numerical experiments with the use of~\eqref{eq:solution:master}, we did not consider estimates for power losses, which can constitute an additional difficulty for convergence.

Additionally, we have tested our methodology without RRP in order to assess its effectiveness in the presented case studies. In the 3-GEN case, all voltage matrices were rank-1. In the IEEE-118 case, on the other hand, the maximum rank of voltage matrices were rank-13 at the final of the iterations.

\begin{table}[h!]
  \centering
  \caption{Convergence of upper and lower bounds for the IEEE-118 case.}
  \pgfplotstableset{every head row/.style={before row=\toprule, after row=\midrule}, every last row/.style={after row=\toprule},
    display columns/0/.style={column name={Iteration}, },
    display columns/1/.style={column name={LB}, },
    display columns/2/.style={column name={UB}, },
    columns/0/.style={fixed, column type=c},
    columns/1/.style={fixed,fixed zerofill,precision=2 ,column type=c},
    columns/2/.style={fixed,fixed zerofill,precision=2,column type=c},
  }
  \begin{filecontents*}{UBLB118bus.txt}
1	1615160.464957 	4884906.437774
2	1615209.077457 	4765981.949478
3	1615191.008669 	3545654.320185
4	1615755.126121 	1700015.570558
5	1699320.683426 	1699482.331184
6	1699482.331184 	1699482.331184

\end{filecontents*}
\pgfplotstabletypeset[columns={ [index]0,[index]1,[index]2}]{UBLB118bus.txt}
    \label{tab:ieee118sol1}%
\end{table}

\section{Conclusion}
\label{sec:conclusion}
In this paper we presented a formulation of the UC problem with AC power flow constraints over the voltage variable space solved by SDP relaxation in a Benders'-like decomposition approach. Our proposed methodology incorporates linearized formulations of constraints of the OPF subproblem into a modified master problem, thus improving the quality of generated cuts. This modification of the UC master problem has allowed convergence of the algorithm in very few iterations, with low rank voltage matrices resulted from our proposed RRP, in both the three-generator small-scale and the IEEE 118-bus-based larger-scale case studies, as illustrated in the numerical experiments. Our Benders'-like decomposition of the problem with convex relaxation of AC OPF constraints enables the incorporation of contingency analysis to account for scenarios of system failure analogously to the work of Fu et al.~\cite{6570740} in future developments. Moreover, because of the nature of the cuts generated from a SDP relaxation, it should be expected such cuts to be more conservative as compared to traditional linear formulations in the literature. Likewise, it is also expected to increase efficiency of the algorithm by adapting strategies proposed in the literature to improve cut selection in linear formulations~\cite{AcelBen,InexactCut,Accelerating,1984327}. Future work also includes an assessment of the feasibility and optimality of the solutions found by our method in comparison to other solvers, as well as the study of the effects of OPF rank-1 solution recovery approaches on the solution to the UC problem.

\appendix

\section{Additional 3-GEN Case Data}

\begin{table}[H]
  \centering
  \caption{General data for the 3-GEN case.}
  \pgfplotstableset{every head row/.style={before row=\toprule, after row=\midrule}, every last row/.style={after row=\bottomrule},columns/T/.style={column name=$T \left[ \si{\hour}\right]$}, columns/Sbase/.style={column name=$Sbase \left[ \si{\mega VA} \right]$}, columns/Vo/.style={column name=$V_{0} \left[ p.u.\right]$}, columns/Vmax/.style={column name=$\overline{V} \left[ p.u.\right]$},columns/Vmin/.style={column name=$\underline{V} \left[ p.u.\right]$},}
  \begin{filecontents*}{3-GENTgldata.dat}
T	Sbase	{$i_{0}$}	Vo	Vmax	Vmin
24	100.000000	1	1.000000	1.050000	0.950000
\end{filecontents*}
\pgfplotstabletypeset[ 1000 sep={,}, columns/4/.style={column name=$\overline{V}$}, columns/info/.style={fixed,fixed zerofill,precision=1,showpos, column type=r,},columns={ [index]0,[index]1,[index]2,[index]3,[index]4,[index]5}]{3-GENTgldata.dat}
\end{table}

\begin{table}[H]
  \centering
  \caption{Transmission line data.}
\pgfplotstabletranspose[string type,
    colnames from=parametro,
    input colnames to=parametro
]\mytablenew{3-GENHlinedata.dat}
 \pgfplotstableset{every head row/.style={before row=\toprule, after row=\midrule}, every last row/.style={after row=\bottomrule},columns/parametro/.style={column name=Line}, columns/Fmax/.style={column name=$\overline{F} \; \left[\si{\mega \watt} \right]$}, columns/r/.style={column name=$r \; \left[p.u.\right]$}, columns/x/.style={column name=$x \; \left[p.u.\right]$},}
\pgfplotstabletypeset[ 1000 sep={,},columns/info/.style={fixed,fixed zerofill,precision=1,showpos,},
    every head row/.style={
        before row=\toprule,
        after row=\midrule
    },
    every last row/.style={
        after row=\bottomrule
    },
string type]{\mytablenew}
\end{table}

\begin{table}[H]
  \centering
  \caption{Load data (in \si{\mega\watt} and \si{\mega}VAR respectively). Load on buses 1,2 and 6 are zero.}
  \pgfplotstableset{ every head row/.style={before row=\toprule, after row=\midrule}, every last row/.style={after row=\bottomrule},}
  \begin{filecontents*}{3-GENTbardata.dat}
tempo	dp1	dq1	dp2	dq2	dp3	dq3	dp4	dq4	dp5	dq5	dp6	dq6
1	0.00	0.00	0.00	0.00	34.05	9.01	68.10	18.01	68.10	18.01	0.00	0.00
2	0.00	0.00	0.00	0.00	31.92	8.45	63.84	16.88	63.84	16.88	0.00	0.00
3	0.00	0.00	0.00	0.00	30.86	8.17	61.71	16.32	61.71	16.32	0.00	0.00
4	0.00	0.00	0.00	0.00	29.80	7.88	59.58	15.76	59.58	15.76	0.00	0.00
5	0.00	0.00	0.00	0.00	29.80	7.88	59.58	15.76	59.58	15.76	0.00	0.00
6	0.00	0.00	0.00	0.00	30.86	8.17	61.71	16.32	61.71	16.32	0.00	0.00
7	0.00	0.00	0.00	0.00	34.05	9.01	68.10	18.01	68.10	18.01	0.00	0.00
8	0.00	0.00	0.00	0.00	40.44	10.70	80.86	21.38	80.86	21.38	0.00	0.00
9	0.00	0.00	0.00	0.00	46.29	12.25	92.57	24.48	92.57	24.48	0.00	0.00
10	0.00	0.00	0.00	0.00	50.54	13.37	101.08	26.73	101.08	26.73	0.00	0.00
11	0.00	0.00	0.00	0.00	52.69	13.94	105.34	27.86	105.34	27.86	0.00	0.00
12	0.00	0.00	0.00	0.00	53.20	14.08	106.40	28.14	106.40	28.14	0.00	0.00
13	0.00	0.00	0.00	0.00	52.69	13.94	105.34	27.86	105.34	27.86	0.00	0.00
14	0.00	0.00	0.00	0.00	53.20	14.08	106.40	28.14	106.40	28.14	0.00	0.00
15	0.00	0.00	0.00	0.00	53.20	14.08	106.40	28.14	106.40	28.14	0.00	0.00
16	0.00	0.00	0.00	0.00	51.61	13.65	103.21	27.35	103.21	27.35	0.00	0.00
17	0.00	0.00	0.00	0.00	51.08	13.52	102.14	27.01	102.14	27.01	0.00	0.00
18	0.00	0.00	0.00	0.00	51.08	13.52	102.14	27.01	102.14	27.01	0.00	0.00
19	0.00	0.00	0.00	0.00	49.48	13.09	98.95	26.17	98.95	26.17	0.00	0.00
20	0.00	0.00	0.00	0.00	48.95	12.95	97.89	25.88	97.89	25.88	0.00	0.00
21	0.00	0.00	0.00	0.00	48.95	12.95	97.89	25.90	97.89	25.90	0.00	0.00
22	0.00	0.00	0.00	0.00	49.48	13.09	98.95	26.17	98.95	26.17	0.00	0.00
23	0.00	0.00	0.00	0.00	46.25	12.24	92.57	24.48	92.57	24.48	0.00	0.00
24	0.00	0.00	0.00	0.00	38.31	10.14	76.61	20.26	76.61	20.26	0.00	0.00
\end{filecontents*}
\pgfplotstabletypeset[ 1000 sep={,},string replace={tempo}{$t$},string replace={Vomin}{$\underline{\nu}$ $\left[\si{\hecto \meter^{3} } \right]$},columns/tempo/.style={column name=$T \left[ \si{\hour}\right]$,string type,column type=l,},columns/info/.style={fixed,fixed zerofill,precision=2,showpos,},columns={[index]0, [index]5, [index]6, [index]7, [index]8, [index]9, [index]10}]{3-GENTbardata.dat}
\end{table}

\section*{Acknowledgements}
This research was funded in part by the São Paulo Research Foundation (FAPESP) under grant numbers 13/07570-7 and 14/08972-4. The authors would like to thank Professor Mohammad Shahidehpour of Illinois Institute of Technology for his invaluable insights in the early versions of this manuscript.




\bibliographystyle{elsarticle-num-names}
\bibliography{arxiv}

\end{document}